\RequirePackage{fix-cm}
\documentclass[a4paper, 12pt, leqno]{article}%{amsart}

\usepackage{latexsym, amsmath, amsthm, amssymb, graphicx, graphics, epsfig, bm}
\usepackage{bm}
\usepackage{tikz-cd}

\usepackage{color}

\usepackage{fullpage}

\usepackage{enumitem}

\newcommand{\comment}[1]{}

\usepackage[T1]{fontenc}
\usepackage{latexsym}
\usepackage{manfnt}
\usepackage[sans]{dsfont}
\usepackage{palatino}
\usepackage[sc, osf]{mathpazo} 
\usepackage{eulervm}
\usepackage{euscript}
\newtheorem{Thm}{Theorem}

\newtheorem{thm}{Theorem}[section]
\newtheorem*{thm*}{Theorem}
\newtheorem*{thmA'}{Theorem A'}
\newtheorem*{thmC'}{Theorem \ref{thm:K1}'}
\newtheorem*{thmD'}{Theorem \ref{thm:inv}'}
\newtheorem{lem}[thm]{Lemma}
\newtheorem*{lem*}{Lemma}

\newtheorem*{cor*}{Corollary}
\newtheorem{prop}[thm]{Proposition}
\newtheorem*{prop*}{Proposition}
\newtheorem*{claim*}{Claim}

\theoremstyle{definition}

\newtheorem*{rem*}{Remark}

\newtheorem{exm}[thm]{Example}

\numberwithin{equation}{section}

\renewcommand{\proof}{\vspace{-8pt}\noindent\textit{\textbf{Proof. }}}

\renewcommand{\endproof}{$\square$}

\newcommand{\txt}[1]{\;\;\text{#1}\;\;}
\newcommand{\q}[1]{{``#1''}}
\newcommand{\br}[1]{\left(#1\right)}
\newcommand{\brr}[1]{\left[#1\right]}
\newcommand{\bra}[1]{\left\langle #1 \right\rangle}

\renewcommand{\leq}{\leqslant}
\renewcommand{\geq}{\geqslant}

\newcommand{\hookto}{\hookrightarrow}
\newcommand{\hookr}{\arrow[hookrightarrow]}

\newcommand{\hookrb}{\arrow[hookrightarrow, color=blue]}

\newcommand{\set}[1]{\left\{#1\right\}}
\newcommand{\sett}[2]{\left\{#1 \,|\, #2\right\}}

\newcommand{\norm}[1]{\left\|#1\right\|}
\newcommand{\case}[1]{\begin{cases}#1\end{cases}}
\newcommand{\inv}{^{-1}}
\newcommand{\restr}[2]{\left. #1 \right|_{#2}}
\newcommand{\eps}{\varepsilon}

\newcommand{\matr}[1]{
\begin{pmatrix}
#1
\end{pmatrix}}

\newcommand{\smatr}[1]{\br{ \begin{smallmatrix}#1\end{smallmatrix} } }

\newcommand{\R}{\mathbb R}

\newcommand{\Z}{\mathbb Z}
\newcommand{\CC}{\mathbb C}

\newcommand{\RP}{\mathbb{R}P}

\DeclareMathOperator{\Res}{Res}
\DeclareMathOperator{\Spec}{Spec}

\DeclareMathOperator{\Ran}{Im}
\DeclareMathOperator{\Ker}{Ker}

\DeclareMathOperator{\Gr}{Gr}

\DeclareMathOperator{\Id}{Id}

\DeclareMathOperator{\dom}{dom}

\newcommand{\sa}{^{\mathrm{sa}}}
\newcommand{\eu}{_{\mathrm{K}}}
\newcommand{\ess}{_{\mathrm{ess}}}

\newcommand{\gla}{^{>\lambda}}

\newcommand{\U}{\mathcal U}

\newcommand{\UK}{\U_K}
\newcommand{\UF}{\U_F}

\newcommand{\sP}{{^s\Proj}}
\newcommand{\sD}{{^s\D}}

\newcommand{\Reg}{\mathcal R}
\newcommand{\Rsa}{\Reg\sa}
\newcommand{\RF}{\Reg_F}
\newcommand{\RK}{\Reg_K}

\newcommand{\rR}{{^r\Reg}}
\newcommand{\gR}{{^g\Reg}}

\newcommand{\taub}{\bar{\tau}}
\newcommand{\sigb}{\bar{\sigma}}

\newcommand{\Rp}{\Reg^+}
\newcommand{\Rm}{\Reg^-}

\newcommand{\Bp}{\B^+}
\newcommand{\Bm}{\B^-}
\newcommand{\RFp}{\Reg_F^+}
\newcommand{\RFm}{\Reg_F^-}
\newcommand{\RKp}{\Reg_K^+}
\newcommand{\RKm}{\Reg_K^-}
\newcommand{\BFp}{\B_F^+}
\newcommand{\BFm}{\B_F^-}

\newcommand{\D}{\mathcal D}
\newcommand{\Deu}{\D\eu}

\newcommand{\inve}{_{\mathrm{inv}}}
\newcommand{\Rinv}{\Reg\inve}
\newcommand{\RKinv}{\Reg_{K,\,\mathrm{inv}}}

\newcommand{\lla}{[-\lambda, \lambda]}

\newcommand{\BH}{\B(H)}

\newcommand{\B}{\mathcal B}
\newcommand{\Bsa}{\B\sa}

\newcommand{\BF}{\B_F}
\newcommand{\BsaH}{\B\sa(H)}

\newcommand{\Proj}{\mathcal P}
\newcommand{\PK}{\Proj_{K}}
\newcommand{\PF}{\Proj_F}

\newcommand{\st}{^{\star}}
\newcommand{\Pst}{\Proj\st}

\newcommand{\Pfin}{\Proj_{\mathrm{fin}}}

\newcommand{\LP}{\Proj^1}

\DeclareMathOperator{\one}{\mathbb{1}}
\newcommand{\onep}{\one_{[0,+\infty)}}

\newcommand{\onem}{\one_{(-\infty,0)}}

\newcommand{\Do}{\mathring{\D}}

\newcommand{\K}{\mathcal K}

\renewcommand{\H}{\mathcal H}

\newcommand{\T}{\mathcal T}

\newcommand{\A}{\mathcal A}

\newcommand{\Ahat}{\hat{A}}
\newcommand{\Hhat}{\hat{H}}

\newcommand{\f}{\bm{\chi}}

\newcommand{\p}{\mathbf p}
\newcommand{\pt}{\tilde{\p}}
\newcommand{\pinf}{p_{\infty}}

\newcommand{\x}{\bar{x}}

\newcommand{\kap}{\bm{\kappa}}
\newcommand{\kat}{\tilde{\kap}}

\newcommand{\uu}{\mathbf u}

\tikzset{cong/.style={draw=none,edge node={node [sloped, allow upside down, auto=false]{$\cong$}}},
         Isom/.style={every to/.append style={edge node={node [above,sloped, inner sep=0.4pt, allow upside down, auto=false]{$\sim$}}}}}

\renewcommand{\thefootnote}{\fnsymbol{footnote}}

\newlength{\otstup} 
\newcommand{\sub}[1]{\vspace{\otstup}\textbf{#1}\hspace*{0.5em}}
\newcommand{\upskip}{\vspace{-\otstup}}

\title{Spaces of unbounded Fredholm operators: \\ I. Homotopy equivalences} 
\author{Marina\,Prokhorova}
\date{}

\begin{document}

\maketitle

\footnotetext{\hspace*{-1.8em}Department of Mathematics, University of Haifa (Israel) \\
Department of Mathematics, Technion --- Israel Institute of Technology \vspace*{0.5em}}
\footnotetext{\hspace*{-1.8em}This work was partially supported by ISF grants no. 431/20, 844/19, and 876/20.}

\renewcommand{\thefootnote}{[\arabic{footnote}]}

\renewcommand{\baselinestretch}{1.03}
\selectfont

\begin{abstract}
This paper is devoted to the space of unbounded Fredholm operators equip\-ped with the graph topology, the subspace of operators with compact resolvents, and their subspaces consisting of self-adjoint operators. Our main results are the following:
(1) Natural maps between these four spaces and classical spaces of bounded operators representing K-theory are homotopy equivalences. This provides an alternative proof of a particular case of results of Joachim.
(2) The subspace of essentially positive Fredholm operators represents odd K-theory.
(3) The subspace of invertible operators in each of these spaces of unbounded operators is contractible.
\end{abstract}

%\let\oldnumberline=\numberline
%\renewcommand{\numberline}{\protect\vspace{-19pt}\oldnumberline}
%\addtocontents{toc}{\protect\renewcommand{\bfseries}{}}

%\tableofcontents

%\section*{Contents}

\bigskip
\textbf{\large{Contents}}
\bigskip

\newcommand{\cont}[3]{\hbox to 0.7\textwidth{#1.\hspace{5pt} #2 \hfil #3}}
\newcommand{\conte}[2]{\hbox to 0.7\textwidth{\hspace{15.5pt} #1 \hfil #2}}

\cont{1}{Introduction}{2}
\cont{2}{Preliminaries: operators, spaces, and maps}{8} 
\cont{3}{Proof of Theorem D\hbox{}: Riesz and norm topology}{12} 
\cont{4}{Proof of Theorem D\hbox{}: graph topology}{14} 
\cont{5}{Proof of Theorem D\hbox{}: odd case}{17} 
\cont{6}{Proof of Theorems B and C}{19} 
\cont{7}{Proof of Theorem A}{23} 
\cont{8}{Essentially positive / negative operators}{25} 
\conte{Conventions and notations}{27} 
\conte{References}{27} 

\bigskip

\section{Introduction}\label{sec:intro}

Let $H$ be a separable complex Hilbert space of infinite dimension.
We denote by $\B(H)$ the space of bounded linear operators on $H$ with the norm topology;
by $\K(H)$, $\U(H)$, and $\Proj(H)$ the subspaces of $\B(H)$ 
consisting of compact operators, unitary operators, and projections respectively 
(by \q{projections} we always mean orthogonal projections, that is, self-adjoint idempotents). 

\sub{Regular operators.}
An unbounded operator $A$ on $H$ is a linear operator defined on a subspace $\dom(A)$ of $H$ and taking values in $H$.
Such an operator $A$ is called closed if its graph is closed in $\Hhat = H\oplus H$, 
and densely defined if its domain $\dom(A)$ is dense in $H$.
It is called \emph{regular} if it is closed and densely defined.
The adjoint $A^*$ of a regular operator is itself a regular operator.

Let $\Reg(H)$ denote the set of all regular operators on $H$
and $\Reg\sa(H)\subset \Reg(H)$ denote the subset of self-adjoint operators.

\sub{Two topologies on regular operators.}
The two most useful topologies on the set $\Reg(H)$ of regular operators are the Riesz topology and the graph topology.

The \emph{Riesz topology} on $\Reg(H)$ is induced by the inclusion $\f\colon\Reg(H)\hookto\B(H)$
from the norm topology on the space $\B(H)$ of bounded operators,
where $\f$ is the so called \q{bounded transform map}, $\f(A) = A(1+A^*A)^{-1/2}$.
The image of this inclusion lies in the closed unit ball $\D(H)$ of $\B(H)$.

The \emph{graph topology} on $\Reg(H)$ is induced by the inclusion 
$\p\colon\Reg(H)\hookto\Proj(\Hhat)$ from the norm topology on the space $\Proj(\Hhat)$ of projections in $\Hhat$, 
where $\p$ is the map taking a regular operator to the orthogonal projection onto its graph.

Bounded operators are regular\footnote{\hspace{1pt}and unbounded, 
but we will refrain from using the term \q{unbounded} when it may lead to confusion.}.
This defines the natural inclusion $\B(H)\hookto\Reg(H)$, 
which is continuous and, moreover, an embedding with respect to both Riesz and graph topology on $\Reg(H)$.

Let $\rR(H)$, resp. $\gR(H)$ denote the space of regular operators on $H$ equipped with the Riesz, resp. graph topology.
The identity map $\rR(H)\to\gR(H)$ is continuous,
so we have a sequence of continuous maps
\begin{equation*}
	\begin{tikzcd}
		\B(H) \hookr{r} & \rR(H) \arrow{r}{\Id} & \gR(H) \hookr{r}{\p} & \Proj(\Hhat).
	\end{tikzcd}
\end{equation*}
Moreover, $\p$ can be written as a composition of $\f$ with a continuous map $\pt\colon\D(H)\to\Proj(\Hhat)$
(see formula \eqref{eq:pt} below), 
making the following diagram commutative:

\begin{equation}
	\begin{tikzcd}
		\B(H) \hookr{r} &
		\rR(H) \arrow{r}{\Id} \hookr{d}{\f} &
		  \gR(H) \hookr{d}{\p} \\
		& \D(H) \arrow{r}{\pt} &
		  \Proj(\Hhat)
	\end{tikzcd}
\end{equation}

\upskip
\sub{Fredholm operators and operators with compact resolvent.}
Similar to the bounded case, a regular operator is called Fredholm 
if its range is closed and its kernel and cokernel are finite-dimensional. 
We denote by $\BF(H)$, resp. $\RF(H)$ the subset of $\B(H)$, resp. $\Reg(H)$ 
consisting of Fredholm operators.

A regular operator is said to have \emph{compact resolvent} 
if both $(1+A^*A)\inv$ and $(1+AA^*)\inv$ are compact operators\footnote{\hspace{1pt}
If the resolvent set of $A$ is non-empty, 
then this definition agrees with the usual sense of the words \q{compact resolvent}:
$(A-\lambda)\inv$ is compact for every $\lambda\in\Res(A)$.
However, a regular operator may have an empty resolvent set; 
the definition we use covers such operators as well.\vspace{4pt}}.
Every operator with compact resolvent is Fredholm.
We denote by $\RK(H)$ the subset of $\Reg(H)$ consisting of regular operators with compact resolvents,
\[ \RK(H) = \sett{A\in\Reg(H)}{(1+A^*A)\inv, (1+AA^*)\inv \in\K(H)}. \]
A self-adjoint regular operator $A$ has compact resolvent if and only if $(A+i)\inv$ is compact;
in this case $A$ has a discrete real spectrum.

\sub{Homotopy type.}
It is a classical result of Atiyah and J\"anich \cite{Atiyah}
that the space $\BF$ of bounded Fredholm operators\footnote{\hspace{1pt}
We will often omit mention of $H$ in order to simplify notations.} 
is a classifying space for the functor $K^0$.
It can be easily seen that in the Riesz topology both $\RF$ and $\RK$ have the same homotopy type as $\BF$.
In contrast with that, the homotopy type of the spaces $\RF$ and $\RK$ 
equipped with the \emph{graph} topology remained unknown for a long time.

In 2003 Joachim \cite{Jo} showed, using results of \cite{BJS}, 
that $\gR_F$ and $\gR_K$ are classifying spaces for the functor $K^0$, similar to the bounded case.
He also proved the $K^1$-analog of this result for regular self-adjoint operators.
Moreover, he proved these results in a more general situation, 
for the Hilbert module $H_A=A\otimes H$ over a unital $C^*$-algebra $A$
(the Hilbert space case corresponds to $A=\CC$).

However, Joachim's proofs are based on a fairly advanced machinery of Kasparov KK-theory, 
even in the case $A=\CC$. 
In this paper we give a transparent proof of this particular case of Joachim's results. 
Moreover, we show that natural maps between these spaces of unbounded operators 
and classical spaces of bounded operators representing $K$-theory 
are homotopy equivalences, the embedding $\B_F\hookto\RF$ being one example. 

Our proof of homotopy equivalence is based on covering of spaces (both the source and the target of a map)
by open subsets in such a way that the \q{gluing pattern} of these coverings is 
exactly the same, up to contractible factors, and is preserved by the map.
The same method allows us to prove that the spaces of essentially positive and essentially negative 
Fredholm operators equipped with the graph topology represent odd K-theory, 
in contrast with the bounded case where they are contractible.

\sub{Spaces and maps.}
A regular operator $A$ is Fredholm if and only if its bounded transform $\f(A)$ is Fredholm; 
$A$ has compact resolvent if and only if $a=\f(A)$ is essentially unitary 
(that is, both $1-a^*a$ and $1-aa^*$ are compact operators).
Let $\D_F(H)$, resp. $\Deu(H)$ denote the subspace of the unit ball $\D(H)$ 
consisting of Fredholm, resp. essentially unitary operators.
Then the bounded transform provides the embeddings 
$\rR_F(H)\hookto \D_F(H)$ and $\rR_K(H)\hookto \Deu(H)$.

Replacing the Riesz topology by the graph topology
and the bounded transform by the embedding $\p\colon\gR(H)\hookto\Proj(\Hhat)$,
one gets a similar picture.
Namely, let $p_0$ and $\pinf$ be the orthogonal projections of $\Hhat=H\oplus H$ onto 
the \q{horizontal} subspace $H\oplus 0$ and the \q{vertical} subspace $0\oplus H$, respectively.
Then $\p$ provides the embeddings
\begin{align*}
 \gR_K(H)\hookto  \PK(\Hhat) &= \sett{p\in\Proj(\Hhat)}{p-\pinf \text{ is compact}}, \\
 \gR_F(H)\hookto  \PF(\Hhat) &= \sett{p\in\Proj(\Hhat)}{p-p_0 \text{ is Fredholm}} 
\end{align*}
(see Proposition \ref{prop:RKF-PKF} and the subsequent comment).
The space $\PF$ is called the Fredholm Grassmanian and $\PK$ is called the restricted Grassmanian.

Taking all this together, we obtain a commutative diagram of continuous maps
\begin{equation}\label{diag:K0}
	\begin{tikzcd}
    & \rR_K(H) \arrow{r}{\Id} \hookrb{d} & \gR_K(H) \hookr{r}{\p} \hookrb{d} & \PK(\Hhat) \hookrb{d}	\\
    \B_F(H) \hookr{r} & \rR_F(H) \arrow{r}{\Id} & \gR_F(H) \hookr{r}{\p} & \PF(\Hhat)	
  \end{tikzcd}
\end{equation}
Our first main result is the following theorem.

\begin{Thm}\label{thm:K0}
All the maps on Diagram \eqref{diag:K0} are homotopy equivalences.
Consequently, all the spaces on the diagram are classifying spaces for the functor $K^0$.
\end{Thm}

\upskip
In fact, we prove slightly stronger result in Section \ref{sec:proofA}, 
adding $\D_K$ and $\D_F$ to the picture.
See Diagram \ref{diag:K0all} and Theorem \ref{thm:K0}'.
We do not include it in the Introduction in order to avoid three-dimensional diagrams here.

\sub{Self-adjoint operators.}
By a classical result of Atiyah and Singer \cite{ASi}, 
the space $\BF\sa$ of bounded self-adjoint Fredholm operators\footnote{\hspace{1pt}
The superscript $\sa$ (resp. subscript $_F$) always stays for \q{self-adjoint} (resp. \q{Fredholm}).}
has three connected components;
two of them, $\BFp$ and $\BFm$, are contractible, 
while the third component $\BF\st$ is a classifying space for the functor $K^1$.
Here $\BFp$ (resp. $\BFm$) is the subspace of $\BF\sa$ 
consisting of essentially positive (resp. essentially negative) operators.
Recall that a bounded self-adjoint operator is called essentially positive (resp. essentially negative)
if it is positive (resp. negative) on some invariant subspace of $H$ of finite codimension.

The definition of essential positivity/negativity for regular self-adjoint operators 
is exactly the same\footnote{\hspace{1pt}
In the bounded case, there is an equivalent definition: a bounded operator is called essentially positive 
if its essential spectrum is contained in the positive ray $[0,+\infty)$.
However, this definition is no longer meaningful for unbounded operators;
for example, the essential spectrum of a regular self-adjoint operator with compact resolvent is always empty.}.
Equivalently, a self-adjoint operator $A$ is essentially positive (resp. essentially negative)
if its bounded transform $\f(A)$ has this property. 
Let $\Rp$ and $\Rm$ denote the corresponding subsets of $\Rsa$, 
and let $\Reg\st = \Rsa\setminus\br{\Rp\cup\Rm}$ denote the subset of operators 
that are neither essentially positive nor essentially negative.

Similarly to the bounded case, the space $\rR\sa_F$ of regular self-adjoint Fredholm operators 
equipped with the \emph{Riesz} topology has three connected components;
two of them, $\rR_F^+$ and $\rR_F^-$, are contractible, 
while the third component $\rR\st_F$ is a classifying space for the functor $K^1$.
The same holds for the subspace $\rR\sa_K$ of operators with compact resolvent: 
it has three components; two of them, $\rR_K^-$ and $\rR_K^+$ are contractible, 
while the third component $\rR\st_K$ is a classifying space for the functor $K^1$.  
The natural embedding $\BF\sa\hookto\rR_F\sa$ is a homotopy equivalence \cite[Theorem 5.10]{Le}
and the bounded transform $\f\colon\rR_F\sa\hookto\Bsa_F$ is homotopy inverse to it.
In particular, both maps preserve partition of the spaces into the three connected components. 

However, with the \emph{graph} topology the situation changes drastically.
The space $^g\Rsa_F$ is path connected \cite[Theorem 1.10]{BLP}, as well as $^g\Rsa_K$.
We show that each of the subspaces $^g\RKm$, $^g\RKp$, and $^g\RK\st$ is dense in $^g\Rsa_K$, 
see Theorem \ref{thm:dens}.

Booss-Bavnbek, Lesch, and Phillips asked in \cite[Remark 1.11]{BLP} 
whether the two \q{trivial parts}, $\RFm$ and $\RFp$, are contractible in the graph topology.
The following theorem, together with Theorem \ref{thm:K1} below, gives a negative answer to this question
and shows that each of these subspaces is a classifying space for the functor $K^1$.

\begin{Thm}\label{thm:ess-pos}
In the following diagrams, all the embeddings are homotopy equivalences in the graph topology:
\begin{equation}\label{diag:+-*}
	\begin{tikzcd}
		\RKm \hookr{r} \hookrb{d} & \RK\sa  \hookrb{d} 
		& &
		\RKp \hookr{r} \hookrb{d} & \RK\sa  \hookrb{d} 
		& &
		\RK\st \hookr{r} \hookrb{d} & \RK\sa  \hookrb{d} 
		\\
		\RFm \hookr{r} & \RF\sa    
		& &
		\RFp \hookr{r} & \RF\sa
		& &
		\RF\st \hookr{r} & \RF\sa
	\end{tikzcd}
\end{equation}
\end{Thm}

We illustrate this result in Section \ref{sec:ess-pos} by an example of a loop of elliptic boundary conditions 
for the differential operator $-d^2/dt^2$ on the interval $[0,1]$; 
the corresponding loop of unbounded self-adjoint operators on $H=L^2[0,1]$ has non-vanishing spectral flow,
although all the operators are essentially positive.

\upskip
\sub{Self-adjoint operators and the Cayley transform.}
The space $\Hhat=H\oplus H$ has a natural (complex) symplectic structure given by the 
symplectic form $\omega(\xi,\eta) = \bra{iI\cdot\xi,\eta}$, 
where $\bra{\cdot,\cdot}$ is the inner product and 
$I = \smatr{0 & -i \\ i & 0}$ is a symmetry (that is, a self-adjoint unitary operator).
We call a projection $p\in\Proj(\Hhat)$ Lagrangian (with respect to $\omega$ or $I$) 
if its range is a Lagrangian subspace of $\Hhat$
(equivalently, the symmetry $2p-1$ anticommutes with $I$, or $Ip = (1-p)I$).
Let $\LP(\Hhat;I)$ denote the subspace of $\Proj(\Hhat)$ consisting of Lagrangian projections.
There is a natural homeomorphism from $\LP(\Hhat;I)$ to the unitary group $\U(H)$; 
we discuss it in Section \ref{sec:maps} in more detail.

A regular operator $A$ is self-adjoint if and only if 
its graph is a Lagrangian subspace in $\Hhat$ with respect to this symplectic structure,
that is, $\p(A)\in\LP(\Hhat;I)$.
This provides a natural embedding $\p\colon\gR\sa(H)\hookto\LP(\Hhat;I)$.
Its composition with the homeomorphism $\LP(\Hhat;I)\cong\U(H)$ 
provides the natural embedding $\kap\colon\gR\sa(H)\hookto\U(H)$,
which is given by the formula $\kap(A) = (A-i)(A+i)\inv$ and is called the Cayley transform.

Restricting the map $\pt\colon\D(H)\to\Proj(\Hhat)$ to the subspace $\D\sa(H)\subset \D(H)$ of self-adjoint operators,
we obtain the continuous map $\kat\colon \D\sa(H)\to\U(H)$,
$\kat(a) = (a-i\sqrt{1-a^2})^2$, which makes the following diagram commutative:
\begin{equation}\label{diag:K1}
	\begin{tikzcd}
		\Bsa \hookr{r} & 
		 \rR\sa \arrow{r}{\Id} \hookr{d}{\f} &
		  \gR\sa \hookr{d}{\kap} \\
		& \D\sa \arrow{r}{\kat} &
		  \U
	\end{tikzcd}
\end{equation}
The Cayley transform provides the embeddings 
\begin{align*}
  \gR_K\sa(H)\hookto\UK(H) &= \sett{u\in\U(H)}{u-1 \text{ is compact}}, \\
  \gR_F\sa(H)\hookto\UF(H) &= \sett{u\in\U(H)}{u+1 \text{ is Fredholm}}.
\end{align*}
Equivalently, in the Grassmanian picture the map $\p$ provides the embedding 
of $\gR_K\sa$ to the restricted Lagrangian Grassmanian $\LP_K=\LP\cap\PK\cong\UK$ 
and the embedding of $\gR_F\sa$ to the Fredholm Lagrangian Grassmanian $\LP_F=\LP\cap\PF\cong\UF$.

Combining the maps discussed above and restricting them
to the homotopically nontrivial connected components of $\B\sa_F$, $\rR\sa_K$, and $\rR\sa_F$,
we obtain the following commutative diagram of continuous maps:
\begin{equation}\label{diag:K1-1}
	\begin{tikzcd}
    & \rR\st_K \arrow{r}\hookrb{d} & \gR\sa_K \hookr{r}{\kap} \hookrb{d} & \UK \hookrb{d}	\\
    \B\st_F \hookr{r} & \rR\st_F \arrow{r} & \gR\sa_F \hookr{r}{\kap} & \UF	
  \end{tikzcd}
\end{equation}

\begin{Thm}\label{thm:K1}
	All the maps on Diagram \eqref{diag:K1-1} are homotopy equivalences.
	Consequently, all the spaces on the diagram are classifying spaces for the functor $K^1$.
\end{Thm}

\upskip
We prove slightly stronger result in Section \ref{sec:proofBC}, 
adding $\Deu\st=\Deu\cap\B\st$ and $\D_F\st=\D_F\cap\B\st$ to the picture.
See Diagram \ref{diag:K1all} and Theorem \ref{thm:K1}'.

\sub{Invertible operators.}
Recall that a regular operator $A$ is called invertible 
if the map $A\colon\dom(A)\to H$ is bijective. 
In this case the inverse map $A\inv\colon H\to\dom(A)\subset H$ is a bounded operator on $H$.
An invertible operator is automatically Fredholm.
The spectrum of a regular operator $A$ is defined as the set of all $\lambda\in\CC$ 
such that $A-\lambda$ is not invertible.

Our proof of Theorems \ref{thm:K0}--\ref{thm:K1} is based on the following key observation.

\begin{Thm}\label{thm:inv}
\
\begin{enumerate}[topsep=-4pt, itemsep=6pt, parsep=3pt, partopsep=0pt, leftmargin=16pt]
	\item Let $X$ be one of the spaces $\rR\st$, $\gR\sa$, $\gR^-$, $\gR^+$, or $\gR\st$,
and let $X_K = X\cap\Reg_K$ be the subspace of $X$ consisting of operators with compact resolvent.
Then the subspaces of $X$ and $X_K$ consisting of invertible operators are contractible.
More generally, the subspaces 
\begin{equation*}
  X\lla = \sett{A\in X}{\Spec(A)\cap\lla = \emptyset} \quad \text{and} \quad X_K\lla = X_K\cap X\lla
\end{equation*}
of $X$ are contractible for every $\lambda\geq 0$.
	\item Let $X=\rR$ or $\gR$. 
Then the subspaces of $X$ and $X_K = X\cap\Reg_K$ consisting of invertible operators are contractible.
More generally, the subspaces 
\begin{equation*}
  X\lla = \sett{A\in X}{\Spec(\Ahat)\cap\lla = \emptyset}, 
	 \quad \text{where} \; \Ahat = \smatr{0 & A^* \\ A & 0}\in\Rsa(\Hhat),	
\end{equation*}
and $X_K\lla = X_K\cap X\lla$ are contractible for every $\lambda\geq 0$.
\end{enumerate}
\end{Thm}

\upskip
We consider various spaces of invertible operators in Sections \ref{sec:Dr}--\ref{sec:D1}
and prove different parts of Theorem \ref{thm:inv} in Propositions 
\ref{prop:rRll}, \ref{prop:gRllK}, \ref{prop:gRll}, \ref{prop:gR1}, and \ref{prop:RBD1}.
Along with it, we prove contractibility of other spaces that we will need in the proofs of Theorems \ref{thm:K0}--\ref{thm:K1}.

\upskip
\sub{Proof of Theorems \ref{thm:K0}--\ref{thm:K1}.} %Theorems \ref{thm:K1} and \ref{thm:ess-pos}.}
Our proof handles all the maps on Diagrams \eqref{diag:K0}, \eqref{diag:+-*}, and \eqref{diag:K1-1} 
at once and gives an alternative proof even for those maps for which other proofs are known.

The proof is based on a theorem of tom Dieck that says, roughly, 
that a map is a homotopy equivalence if it is locally a homotopy equivalence. 
We need a particular case of \cite[Theorem 1]{tD}, which is stated as follows.
Let $\varphi\colon X\to Y$ be a continuous map.
Let $(X_{\tau})$, resp. $(Y_{\tau})$, be a numerable covering of $X$, resp. $Y$,
indexed by the same index set $\T$.  
Assume that $\varphi(X_{\tau})\subset Y_{\tau}$ and that for every finite subset $T\subset\T$ 
the restriction map $\bigcap_{\tau\in T} X_{\tau}\to \bigcap_{\tau\in T}Y_{\tau}$
is a homotopy equivalence. 
Then $\varphi$ itself is a homotopy equivalence.

Let us describe an idea of the proof on the example of the inclusion map 
\[ \varphi\colon\rR\st_K(H)\to\gR\sa_K(H). \]
Theorem \ref{thm:inv} suggests to choose coverings of $X=\rR\st_K$ and $Y=\gR\sa_K$ 
by contractible open subspaces indexed by real numbers $\lambda$
and consisting of operators $A$ such that $A-\lambda$ is invertible.
We prefer to use finer coverings $(X_{\tau})$ and $(Y_{\tau})$
which are closed under finite intersections and are suitable for the proof of Theorem \ref{thm:K0} as well.
We index them by finite symmetric
(with respect to zero) non-empty subsets $\tau = \set{\tau_1<\ldots<\tau_n}$ of $\R$,  
with $X_{\tau}$, resp. $Y_{\tau}$ being the subspace of $X$, resp. $Y$ 
consisting of operators whose resolvent set contains $\tau$.
We need only to show that the restriction 
$\varphi_{\tau}\colon X_{\tau}\to Y_{\tau}$ of $\varphi$ is a homotopy equivalence for every such $\tau$.

With every operator $A\in X_{\tau}$ we associate 
its \q{finite part} $A' = \restr{A}{V}$ and the \q{infinite part} $A'' = \restr{A}{V^{\bot}}$,
where $V$ is the finite-dimensional range of the spectral projection $\one_{\taub}(A)$ 
corresponding to the convex hull $\taub = [\tau_1,\tau_n]$ of $\tau$
(here $\one_S$ denotes the characteristic function of a subset $S\subset\CC$).
In such a way, we provide $X_{\tau}$ with the natural structure of 
a fiber bundle over the base space $X'_{\tau}$ consisting of pairs $(V,A')$, 
where $V$ is a finite-dimensional subspace of $H$ 
and $A'\in\Bsa(V)$ has the spectrum contained in the finite union 
$\taub\setminus\tau = \cup_{i=1}^{n-1}(\tau_i,\tau_{i+1})$ of open intervals.
The fiber of $X_{\tau}$ over $(V,A')$ is the subspace of $\rR\st_K(V^{\bot})$ 
consisting of operators $A''$ with $\Spec(A'')\cap\taub = \emptyset$.

The space $Y_{\tau}$ is equipped with a fiber bundle structure over the same base space in exactly the same manner;
the only difference is that $\rR\st_K(V^{\bot})$ is replaced by $\gR\sa_K(V^{\bot})$ in the description of the fibers. 

The fiber bundles $X_{\tau}\to X'_{\tau}$ and $Y_{\tau}\to X'_{\tau}$ are locally trivial; 
by Theorem \ref{thm:inv} their fibers are contractible.
It follows that $\varphi_{\tau}$ is a homotopy equivalence.

All the other spaces on our diagrams are handled similarly, 
with an appropriate choice of coverings.

\sub{Acknowledgments.}
This work was done during my postdoctoral fellowship at the Technion -- Israel Institute of Technology; its revised version was prepared during my postdoctoral fellowship at the University of Haifa.

I am grateful to N.V.~Ivanov for many helpful discussions and remarks.
I am also grateful to the anonymous reviewer for careful reading of my paper 
and suggestions helping to improve the exposition.

\section{Preliminaries: operators, spaces, and maps}\label{sec:maps}

In this section we recall basic notions and facts about regular operators
and prove several simple facts that we will use further.
In addition, we explain how the Cayley transform $\kap\colon\Rsa(H)\hookto\U(H)$ 
is related to the projection map $\p\colon\Reg(H)\hookto\Proj(\Hhat)$.

\sub{Adjoint operators.}
The \emph{adjoint operator} of a regular operator $A$ 
is an unbounded operator $A^*$ on $H$ with the domain 
\[ \dom(A^*) = \sett{x\in H}{\text{ there exists } y\in H \text{ such that } 
   \bra{Az,x}=\bra{z,y} \text{ for all } z\in \dom(A)}. \]
For $x\in\dom(A^*)$ such an element $y$ is unique and $A^*x=y$ by definition.
The adjoint of a regular operator is itself a regular operator.
An operator $A\in\Reg(H)$ is called \emph{self-adjoint} if $A^*=A$ (in particular, $\dom(A^*) = \dom(A)$).

\sub{Bounded transform.}
Let us recall the main properties of the bounded transform.
See, e.g., \cite[Theorem 7.5]{Schm} for the detail.

For a regular operator $A$, the operator $1+A^*A$ is regular, self-adjoint, and surjective;
its inverse $(1+A^*A)\inv$ is a bounded self-adjoint operator. % \cite[Theorem V.3.24]{Kato}.
The \emph{bounded transform} (or the Riesz map) 
\[ \f\colon\Reg(H)\to\B(H), \quad \f(A) = A(1+A^*A)^{-1/2}, \] 
defines the inclusion of the set $\Reg(H)$ of regular operators to the closed unit ball 
\[ \D(H)=\sett{a\in\B(H)}{\norm{a}\leq 1} \]
in the space $\B(H)$ of bounded operators.
The image of this inclusion, 
\begin{equation}\label{eq:Do1}
	\Do(H) = \f(\Reg(H)) = \sett{ a\in \D(H)}{1-a^*a \text{ is injective} }, 
\end{equation}
is dense in $\D(H)$.
The inverse map $\f\inv\colon\Do(H)\to\Reg(H)$ is given by the formula 
$\f\inv(a) = a(1-a^*a)^{-1/2} = (1-aa^*)^{-1/2}a$.

If a regular operator $A$ is self-adjoint, then so is $\f(A)$; more generally, $\f(A^*) = \f(A)^*$.
If $a=\f(A)$, then $(1+A^*A)\inv = 1-a^*a$.
In particular, $A$ has compact resolvent if and only if $\f(A)$ is essentially unitary.

\begin{prop}\label{prop:conv}
  The image $\Do(H) = \f(\Reg(H))$ is a convex subset of $\D(H)$.
\end{prop}

\proof
Since $1-a^*a$ is positive for all $a\in\D$, \eqref{eq:Do1} can be written equivalently as
\begin{equation*}\label{eq:Do2}
 \Do = \set{ a\in\D\colon \norm{a\xi} < \norm{\xi} \text{ for every non-zero } \xi\in H}.
\end{equation*}
Let $a_0,a_1\in\Do$ and $s\in[0,1]$. 
Then $a_s=(1-s)a_0+sa_1$ satisfies the inequality
\[
 \norm{a_s \xi} \leq (1-s)\norm{a_0 \xi} + s \norm{a_1 \xi} < \norm{\xi} 
\]
for every non-zero $\xi\in H$, and thus $a_s\in\Do$ as well.
\endproof

\begin{prop}\label{prop:conv+}
For $\lambda\geq 0$, let $\Reg\gla\subset\Reg\sa$ denote the set of self-adjoint operators in $H$
	whose spectrum is contained in the open ray $(\lambda, \infty)$. 
	Then the images under $\f$ of both $\Reg\gla$ and $\RK\gla=\RK\cap\Reg\gla$ are convex subsets of $\D\sa$,
	so $\rR\gla$ and $\rR_K\gla$ are contractible.
\end{prop}

\proof
The first image $\f(\Reg\gla)$ is the intersection of the convex set $\Do$ with %the set
\[ \B^{>\mu} = \sett{a\in\B\sa}{\Spec(a)\subset(\mu,\infty)}, \txt{where} \mu=\f(\lambda). \]
%Let $\mu=\f(\lambda)$ and $\B^{>\mu} = \sett{a\in\B\sa}{\Spec(a)\subset(\mu,\infty)}$.
If $a_0,a_1\in\B^{>\mu}$, then $a_i\geq\mu_i$ for some $\mu_0,\mu_1>\mu$, 
so every convex combination of $a_0$ and $a_1$ is $\geq\min\set{\mu_0,\mu_1}>\mu$ and thus lies in $\B^{>\mu}$.
Hence both $\B^{>\mu}$ and $\f(\Reg\gla)$ are convex.
The second image $\f(\RK\gla)$ is the intersection of two convex sets, $\f(\Reg\gla)$ and 
\[ \sett{a\in\D\sa(H)}{1-a^2\in\K(H) \text{ and } a\geq 0} = 
   \sett{a\in\D\sa(H)}{1-a\in\K(H) \text{ and } a\geq 0}, \]
so it is also convex.	
By the definition of the Riesz topology, $\rR\gla$ and $\rR_K\gla$ are homeomorphic to their bounded transforms,
which are convex and thus contractible.
\endproof

\begin{rem*}
	The same reasoning works for every $\lambda\in\R$, and also for closed rays $[\lambda, \infty)$. 
	But the statement of Proposition \ref{prop:conv+} is enough for our purposes.
\end{rem*}

\sub{Projections.}
The orthogonal projection onto the graph of a regular operator $A$ is given by the formula
\begin{equation}\label{eq:pA}
 \p(A) = \matr{(1+A^*A)^{-1} & (1+A^*A)^{-1}A^* \\ A(1+A^*A)^{-1} & 1-(1+AA^*)^{-1}} \in\Proj(\Hhat). 
\end{equation}
Recall that we denoted by $\rR(H)$ and $\gR(H)$ the space of regular operators on $H$
equipped with the Riesz and the graph topology respectively.
The identity map $\rR(H)\to\gR(H)$ is continuous, that is, $\p\colon\rR(H)\to\Proj(\Hhat)$ 
factors through a continuous map $\pt\colon\Do(H)\to\Proj(\Hhat)$, $\p=\pt\circ\f$.
Moreover, $\pt$ can be (uniquely) continuously extended to the whole closed unit ball $\D(H)$,
making the following square commutative:
\begin{equation*}
	\begin{tikzcd}
		\rR \arrow{r}{\Id} \hookr{d}{\f} &
		  \gR \hookr{d}{\p} \\
		\D \arrow{r}{\pt} &
		  \Proj 
	\end{tikzcd}
\end{equation*}
Such an extension $\pt\colon\D(H)\to\Proj(\Hhat)$ is given by the formula\footnote{\hspace{1pt}
The equality $\sqrt{1-a^*a}\,a^* = a^*\sqrt{1-aa^*}$ follows 
from approximation of the square root on the interval $[0,1]$ by polynomials 
and the fact that $f(1-a^*a)\,a^* = a^*f(1-aa^*)$ for every polynomial $f$.} 
\begin{equation}\label{eq:pt}
 \pt(a) = \matr{1-a^*a & \sqrt{1-a^*a}\,a^* \\ a\sqrt{1-a^*a} & aa^*}
        = \matr{1-a^*a & a^*\sqrt{1-aa^*} \\ a\sqrt{1-a^*a} & aa^*}. 
\end{equation}
We will use two special projections in $\Hhat = H\oplus H$,
onto the \q{horizontal} subspace $H\oplus 0$ and the \q{vertical} subspace $0\oplus H$:
\[ p_0 = \matr{1 & 0 \\ 0 & 0} \text{ and } \pinf = \matr{0 & 0 \\ 0 & 1}. \]

\begin{prop}\label{prop:RKF-PKF}
For every $a\in\D(H)$ the following hold:
\begin{enumerate}[topsep=-4pt, itemsep=0pt, parsep=3pt, partopsep=0pt]
	\item $a$ is essentially unitary if and only if $\pt(a)-\pinf$ is compact.
	\item $a$ is Fredholm if and only if $\pt(a)-p_0$ is Fredholm.
\end{enumerate}
\end{prop}

\upskip
It follows that $A\in\Reg(H)$ has compact resolvent if and only if $\p(A)-\pinf$ is compact;
$A$ is Fredholm if and only if $\p(A)-p_0$ is Fredholm. 

\medskip
\proof
The first part is obvious from formula \eqref{eq:pt}. 
The second part follows from the factorization
\begin{equation*}\label{eq:pta-p0}
	\pt(a)-p_0 = \matr{-a^*a & a^*\sqrt{1-aa^*} \\ a\sqrt{1-a^*a} & aa^*} 
	  = \matr{-a^* & 0 \\ 0 & a}  \cdot \matr{a & -\sqrt{1-aa^*} \\ \sqrt{1-a^*a} & a^*}
\end{equation*}
and the fact that the second factor is unitary and thus invertible.
\endproof

\sub{Cayley transform.}
The Cayley transform 
\[ \kap\colon\Rsa(H)\to\U(H), \quad \kap(A) = (A-i)(A+i)\inv \] 
is a continuous embedding\footnote{\hspace{1pt}See equality \eqref{eq:p-kap} below.} 
of $\gR\sa(H)$ into the unitary group $\U(H)$. 
For self-adjoint $A$ and $a=\f(A)$, the equality $(1+A^2)\inv = 1-a^2$ implies
\[ \kap(A) = \frac{A-i}{A+i} = \frac{a-i\sqrt{1-a^2}}{a+i\sqrt{1-a^2}} = \br{a-i\sqrt{1-a^2}}^2 = \kat(a), \]
where $\kat\colon[-1,1]\to\U(\CC) = \sett{z\in\CC}{|z|=1}$ is a continuous function 
given by the formula 
\begin{equation}\label{eq:kat}
	\kat(a) = \br{a-i\sqrt{1-a^2}}^2 = 2a^2-1-2ia\sqrt{1-a^2}.
\end{equation}
Hence the Cayley transform factors through the bounded transform: $\kap = \kat\circ\f$.
Moreover, $\kat$ can be (uniquely) continuously extended to the whole $\D\sa(H)$;
the corresponding map is given by the same formula \eqref{eq:kat}.

An operator $a\in\D\sa(H)$ is Fredholm if and only if 
\[ \kat(a)\in\UF(H) = \sett{u\in\U(H)}{u+1 \text{ is Fredholm}}; \]
it is essentially unitary (that is, $1-a^2\in\K(H)$) if and only if 
\[ \kat(a)\in\U_K(H) = \sett{u\in\U(H)}{u-1 \text{ is compact}}. \]
This can be seen from the identities
\[ \kat(a)-1 = -2i\sqrt{1-a^2}(a-i\sqrt{1-a^2}) \txt{ and } \kat(a)+1 = 2a\br{a-i\sqrt{1-a^2}} \]
and the fact that $a-i\sqrt{1-a^2}\in\U(H)$.

\sub{Lagrangian projections.}
As we discussed in the Introduction,
the space $\Hhat$ has a natural symplectic structure given by the symmetry $I = \smatr{0 & -i \\ i & 0}$.
A regular operator $A$ is self-adjoint if and only if 
its graph is a Lagrangian subspace in $\Hhat$ with respect to $I$, that is, 
\[ \p(A)\in\LP(\Hhat;I) = \sett{p\in\Proj(\Hhat)}{I(2p-1)+(2p-1)I = 0}. \] 
If $a\in\D(H)$ is self-adjoint, then $\pt(a)$ is Lagrangian with respect to $I$.
But the converse is no longer true, in contrast with regular operators.
For example, $\pt(a)=\smatr{0 & 0 \\ 0 & 1}\in\LP$ for every unitary operator $a$.

Let us consider the grading symmetry $J=\smatr{1 & 0 \\ 0 & -1}$ of $\Hhat$.
A symmetry $r$ anticommutes with $J$ if and only if it has the form 
$r=\smatr{0 & u^* \\ u & 0}$ for some $u\in\U(H)$.
The map $r\mapsto u$ determines a natural homeomorphism $\LP(\Hhat;J)\to\U(H)$.

Take a unitary $v\in\U(\Hhat)$ such that $J=vIv^*$.
Then the conjugation by $v$ takes $\LP(\Hhat;I)$ to $\LP(\Hhat;J)$ 
and thus determines a homeomorphism 
\[ \psi_v\colon\LP(\Hhat;I)\to\U(H). \]
If both $v,v'\in\U(\Hhat)$ conjugate $I$ with $J$, 
then $v' = \smatr{w & 0 \\ 0 & w'}\cdot v$ for some unitaries $w,w'\in\U(H)$.
It follows that 
\begin{equation}\label{eq:psiv}
	\psi_{v'}(p) = w'\cdot\psi_v(p)\cdot w^* \;\text{ for every }\; p\in\LP(\Hhat;I).
\end{equation}
Conversely, for a fixed $v$, every pair $w,w'\in\U(H)$ gives rise to $v'$ satisfying \eqref{eq:psiv}.
Since the unitary group is path connected, 
the isotopy class of the homeomorphism $\psi_v$ does not depend on the choice of $v$.

It is convenient to take $v = \frac{1}{\sqrt{2}}\smatr{-1 & i \\ 1 & i}$.
Then the composition
\begin{equation*}
	\begin{tikzcd}
		\D\sa(H) \arrow{r}{\pt} & \LP(\Hhat;I) \arrow{r}{\psi_v} & \U(H)
 \end{tikzcd}
\end{equation*}
coincides with $\kat$.
In particular, 
\begin{equation}\label{eq:p-kap}
	\psi_v(\p(A)) = \kap(A) \txt{for every} A\in\Rsa(H).
\end{equation}
The homeomorphism $\psi_v$ takes $\pinf$ to $1\in\U(H)$ and $p_0$ to $-1\in\U(H)$
and thus provides homeomorphisms 
\begin{equation}\label{eq:LPKUK}
	 \LP_K = \LP(\Hhat;I)\cap\PK(\Hhat)\to\UK(H) \quad\text{and}\quad 
	   \LP_F = \LP(\Hhat;I)\cap\PF(\Hhat)\to\UF(H). 
\end{equation}
The space $\LP$ is called the Lagrangian Grassmanian, $\LP_F$ the Fredholm Lagrangian Grassmanian,
and $\LP_K$ the restricted Lagrangian Grassmanian.

Composition of $\pt$ with homeomorphisms \eqref{eq:LPKUK} %and taking Proposition \ref{prop:RKF-PKF} into account, 
provides the maps $\kat\colon\Deu\sa(H)\to\UK(H)$ and $\kat\colon\D\sa_F(H)\to\UF(H)$ discussed above.

\section{Proof of Theorem \ref{thm:inv}: Riesz and norm topology}\label{sec:Dr}

In this section we prove the first part of Theorem \ref{thm:inv} for the case $X=\rR\st$, 
as well as its analogue for bounded operators which we will need in Section \ref{sec:proofBC}.

\sub{Canonical decomposition.} 
Recall a standard construction that we will use throughout the paper.
The trivial Hilbert bundle over $\Proj(H)$ with the fiber $H$ 
is canonically decomposed into the orthogonal sum 
\begin{equation}\label{eq:HHH}
	\Proj(H)\times H = \H'\oplus\H'' 
\end{equation}
of two Hilbert bundles, whose fibers are $\H'_p = \Ran p$ and $\H''_p = \Ker p$.
This decomposition is locally trivial in the following sense: 
for every $q\in\Proj(H)$ and the open ball $W$ of radius $1$ around $q$,
there is a continuous map 
\begin{equation}\label{eq:gWU}
	 g\colon W\to\U(H) \;\text{ such that }\; g(p)\cdot p\cdot g(p)^*\equiv q \;\text{ for }\; p\in W.
\end{equation}
See \cite[Proposition 5.2.6]{WO}.

Let $\Pst(H)$ be the subspace of $\Proj(H)$ consisting of projections of infinite rank and corank.
The restrictions of $\H'$ and $\H''$ to $\Pst(H)$ are locally trivial bundles over a paracompact space, 
and their structure group is the unitary group $\U$ with the norm topology, 
which is contractible by Kuiper's theorem \cite{Kui}. 
Thus these restrictions are trivial as Hilbert bundles with the structure group $\U$.
This can be stated as follows: there is a map 
\begin{equation}\label{eq:gP*U}
	 g\colon \Pst(H)\to\U(H) \;\text{ such that }\; g(p)\cdot p\cdot g(p)^*\equiv q 
	\txt{for} p\in\Pst(H),
\end{equation}
where $q\in\Pst(H)$ is some fixed projection.

\sub{Riesz topology.}
The following proposition proves the first part of Theorem \ref{thm:inv} for the case $X=\rR\st$.

\begin{prop}\label{prop:rRll}
 The spaces $\rR\st\lla$ and $\rR\st_K\lla$ are contractible for every $\lambda\geq 0$.
\end{prop}

\proof
The positive spectral projection determines a continuous map 
\[ \onep\colon\rR\st\lla\to\Pst(H), \]
as can be seen from the equality $\onep(A)=\onep(\f(A))$.
The maps 
\begin{equation}\label{eq:onep-rR}
	\rR\st\lla\to \Pst \txt{and} \rR\st_K\lla\to \Pst
\end{equation}
are trivial fiber bundles; their trivializations are given by the maps
\begin{align*} 
 & \rR\st\lla\to F\times\Pst \txt{and} \rR\st_K\lla\to F_K\times\Pst, \\
 & A\mapsto (g(p)\cdot A\cdot g(p)^*,\, p), \quad p=\onep(A),
\end{align*} 
where $g\colon\Pst\to\U$ is a map satisfying \eqref{eq:gP*U}
and $F$, $F_K$ are fibers of the maps \eqref{eq:onep-rR} over $q\in\Pst$.
Let $H' = \Ran q$ and $H'' = \Ker q$.
The map taking $A\in F$ to the pair $\br{\restr{A}{H'}, -\restr{A}{H''}}$ 
determines homeomorphisms 
\[ F\cong \rR\gla(H')\times\rR\gla(H'') \txt{and} F_K\cong \rR_K\gla(H')\times\rR_K\gla(H''), \]
where $\Reg\gla$ denotes the subset of operators whose spectrum is contained in the open ray $(\lambda,+\infty)$.
By Proposition \ref{prop:conv+}, each of the factors is contractible,
so $F$ and $F_K$ are contractible.
The base space $\Pst(H)$ is contractible as well, see \cite[proof of Lemma 3.6]{ASi}.
Therefore, the total spaces $\rR\st\lla$ and $\rR\st_K\lla$ are also contractible.
\endproof

\sub{Norm topology.}
In the proof of Theorems \ref{thm:ess-pos} and \ref{thm:K1}
we will need the following analogue of Theorem \ref{thm:inv} for bounded operators.

\begin{prop}\label{prop:BDUll}
 The following spaces are contractible in the norm topology for every $\lambda\geq 0$:
\begin{align*}
  \B\st\lla &= \sett{A\in \B\st}{\Spec(A)\cap\lla = \emptyset}, \\
	\D\st\lla &= \sett{A\in \D\st=\D\cap\B\st}{\Spec(A)\cap\f\br{\lla} = \emptyset}, \\
  \U\lla &= \sett{A\in \U}{\Spec(A)\cap\kap\br{\lla} = \emptyset}, \\
	\Deu\st\lla &= \Deu\cap\D\st\lla, \\
	\UK\lla &= \UK\cap\U\lla.
\end{align*}
\end{prop}

\proof
For $\B\st\lla$, $\D\st\lla$, and $\D\st_K\lla$ the contractibility is proved 
in exactly the same manner as in Proposition \ref{prop:rRll}, 
with the only difference that we have a convex structure on the fibers from the start, 
without applying the bounded transform.

A contraction of $\U\lla$ to $\set{1}$ is given by the homotopy
\[ h\colon\U\lla\times[0,1]\to\U\lla, \quad h_t(u) = \exp(t\log(u)), \] 
where $\log\colon\set{z\in\CC\colon|z|=1, z\neq -1} \to (-i\pi,i\pi)\subset i\R$ 
is the branch of the natural logarithm.
It preserves $\U_K\lla$ and thus defines a contraction of $\U_K\lla$ to $\set{1}$ as well.
\endproof

\section{Proof of Theorem \ref{thm:inv}: graph topology}\label{sec:Dg}

This section is devoted to the proof of the first part of Theorem \ref{thm:inv} for spaces $X$ from the following list:
\begin{equation}\label{eq:XgR}	
  \gR\sa, \; \gR^-, \; \gR^+, \; \gR\st.
\end{equation}
We complete this task in the end of the section, see Propositions \ref{prop:gRllK} and \ref{prop:gRll}.

\sub{Contraction of $\K(H)$.}
We start from contracting $\K(H)$. 
Of course, this can be done easily by a linear deformation, but our more elaborated construction has an important advantage:
it preserves every subset of $\K(H)$ provided this subset is invariant 
under conjugation by unitary operators and multiplication by numbers from the interval $(0,1]$.
We then apply this deformation to contracting the space $X_K$ for each $X$ from the list above.

We use the approach of Dixmier and Douady to the proof of \cite[Lemma 3]{DD}, 
but modify it in a way which provides a contraction of $\K(H)$.
Let us identify $H$ with $L^2[0,1]$, and let $p_t$ and $q_t$, $t\in[0,1]$, 
be the projections onto the subspaces $L^2[0,t]$ and $L^2[t,1]$ of $H$, respectively.
Then $p_t+q_t\equiv 1$.
These projections can be written as follows: 
\[ p_t=u_t u_t^* \text{ for } t\in(0,1] \txt{and} q_t=v_t v_t^* \text{ for } t\in[0,1), \] 
where $(u_t)_{t\in(0,1]}$ and $(v_t)_{t\in[0,1)}$ are two families of isometries of $H$,
\begin{equation}\label{eq:utvt}
	u_t(f)(s) = \case{ \frac{1}{\sqrt{t}}f(\frac{s}{t}) & \text{ for } s\leq t \\ 0 & \text{ for } s>t}
		\quad\txt{and}\quad
		v_t(f)(s) = \case{ 0 & \text{ for } s\leq t \\ \frac{1}{\sqrt{1-t}}f(\frac{s-t}{1-t}) & \text{ for } s>t} \
\end{equation}
(in terms of \cite[Lemma 2]{DD}, $u_t^* = U_t P_t$).
The operators $u_t$, $u_t^*$, $v_t$, and $v_t^*$ 
depend continuously on $t$ in the strong operator topology and $u_1=v_0=1$.

\begin{prop}\label{prop:K-contr}
	For every $\beta\in\K(H)$, the formula
\begin{equation}\label{eq:ZKtoB}
	h_0(\alpha)=\alpha, \quad h_1(\alpha)=\beta, \quad 
	 h_t(\alpha) = tu_t \beta u_t^* + (1-t)v_t \alpha v_t^* \;\text{ for }\; 0<t<1
\end{equation}
determines a contraction $h\colon[0,1]\times\K(H)\to\K(H)$ of $\K(H)$ to $\set{\beta}$.
\end{prop}

\proof
Since $\alpha$, $\beta$ are compact and $u_t$, $v_t$ are uniformly bounded, the maps 
\[ h'\colon (0,1]\to\K(H),\; h'_t = tu_t \beta u_t^* \txt{and} 
 h''\colon [0,1)\times \K(H)\to \K(H),\; h''_t(\alpha) = (1-t)v_t \alpha v_t^* \]
are continuous\footnote{\hspace{1pt}
Both $u$ and $v$ are maps from an interval to $\sD$, 
the unit ball in $\B(H)$ equipped with the strong operator topology.
Since every compact operator is a limit of finite rank operators,
the multiplication map $\sD\times\K(H)\to\K(H)$, $(a,k)\mapsto ak$ is continuous.
Passing to the adjoint $(ak^*)^*=ka^*$ and multiplying once more from the left, 
we see that the map $\sD\times\K\times\sD\to\K$, $(a,k,b)\mapsto bka^*$ is continuous.}.
In addition, 
\[ \norm{h''_t(\alpha)}\leq (1-t)\norm{\alpha}\to 0 \txt{as} t\to 1, \]
and similarly $\norm{h'_t}\leq t\norm{\beta}\to 0$ as $t\to 0$.
Therefore, $h$ is continuous on the whole domain 
and thus provides a contraction of $\K(H)$ to $\set{\beta}$.
\endproof

\sub{Invertible self-adjoint operators.}
The goal of this subsection is to prove contractibility in the graph topology 
of the space $\Rinv\sa$ of invertible self-adjoint regular operators,
as well as its subspace $\RKinv\sa$ consisting of operators with compact resolvent.

\begin{prop}\label{prop:inverse}
	The map 	
	\begin{equation}\label{eq:inv}
		\gR\inve(H)\to\BH, \quad A\mapsto A\inv 
	\end{equation}
	is continuous. It provides homeomorphisms 
	\[ ^g\Rinv\sa(H)\to Z = \sett{\alpha\in\BsaH}{\alpha \text{ is injective}} \;\text{ and }\; 
	    ^g\RKinv\sa(H)\to Z_K = Z\cap\K(H), \]
	and takes $\Rsa\lla$ onto $\sett{\alpha\in Z}{\norm{\alpha}<\lambda\inv}$.
\end{prop}

\proof 
Clearly, $A\inv$ is injective. 
For every $\alpha\in Z$, the range of $\alpha$ is dense and the graph of $\alpha$ is closed, 
so $\alpha\inv\in\Rinv\sa$.
The bounded linear transformation $x\oplus y \mapsto y\oplus x$ of $H\oplus H$ 
takes the graph of $A$ to the graph of $A\inv$, so \eqref{eq:inv} is graph continuous.
Since the graph topology on $\BH$ coincides with the norm topology, 
this map provides a homeomorphism $\gR\inve\sa\to Z$.
By the definition, an operator $A\in\Rinv\sa$ has compact resolvent if 
$(1+A^2)\inv = \alpha^2(\alpha^2+1)\inv$ is compact;
since $\alpha^2+1$ is an isomorphism, this is equivalent to compactness of $\alpha$. %$a\in\K(H)$
The last statement of the proposition is obvious.
\endproof

\begin{prop}\label{prop:RKinv}
	The space $Z_K$ is contractible, so $^g\RKinv\sa$ is also contractible. % in the graph topology.
\end{prop}

\proof
Choose $\beta\in Z_K$ and let $h$ be a homotopy from Proposition \ref{prop:K-contr}.
The sum $h_t(\alpha) = h'_t+h''_t(\alpha)$ in \eqref{eq:ZKtoB} 
is the orthogonal sum corresponding to the decomposition 
$H = \Ran p_t\oplus \Ran q_t$. 
For $\alpha$ injective, the summands $h'_t$ and $h''_t(\alpha)$
are injective as operators on $\Ran p_t$ and $\Ran q_t$, respectively.
Therefore, $h_t$ preserves $Z_K$ and the corresponding restriction of $h$ 
is a contraction of $Z_K$ to $\set{\beta}$.
The second part of the proposition follows from the first part and Proposition \ref{prop:inverse}.
\endproof

The corresponding contraction of $^g\RKinv\sa$ 
to an arbitrary point $B\in\RKinv\sa$ is given by the formula
\begin{equation}\label{eq:RKtoB}
	H_0(A)=A, \quad H_1(A)=B, \quad H_t(A) = t\inv u_t Bu_t^* + (1-t)\inv v_t Av_t^* \;\text{ for }\; 0<t<1.
\end{equation}

\begin{prop}\label{prop:Rinv}
	The space $Z$ is contractible, so $^g\Rinv\sa$ is also contractible. % in the graph topology.
\end{prop}

\upskip
Cf. \cite[Proposition 5.8]{Le}, 
where the path connectedness of $^g\Rinv\sa$ is shown in a completely different manner.
\medskip

\proof
Let $k\in\K(H)$ be a positive injective compact operator. 
Then 
\begin{equation}\label{eq:bt} 
 c_t=(1-t)+tk 
\end{equation}
is also a positive injective operator for every $t\in[0,1]$. 
The map 
\[ h'\colon [0,1]\times Z\to Z, \quad h'_t(\alpha) = c_t \alpha c_t \] 
is norm continuous, $h'_0=\Id$, and $h'_1(Z)\subset Z_K$.
Since $Z_K$ is contractible, this implies contractibility of $Z$. 
The second part of the proposition follows from the first part and Proposition \ref{prop:inverse}.
\endproof

The corresponding deformation $H'$ of $\gR\inve\sa$ to a subspace of $^g\RKinv\sa$ 
is given by the formula
\begin{equation}\label{eq:RtoRK} 
	H'_t(A) = (c_t A\inv c_t)\inv % \;\text{where}\; a=A\inv. %C_t = (c_t)\inv.
\end{equation}
(the bounded self-adjoint operator $c_t A\inv c_t$ is injective, so its inverse is a regular operator).

\begin{rem*}
Another possible approach to proving the contractibility of spaces $^g\Rinv\sa$ and $^g\RKinv\sa$ 
is to use their natural map $\onep$ to the space $\sP$ of projections with the strong topology. 
The space $\sP$ is contractible, as one can prove using the methods of \cite{DD}, 
and the fibers are convex and hence contractible. 
Unfortunately, these maps seem to fail to be locally trivial bundles, so additional arguments are required. 
The author hopes to return to this topic on another occasion.
Our current approach to the proof, 
namely the direct application of methods of Dixmier and Douady \cite{DD} to $^g\RKinv\sa$ instead of $\sP$, 
allows to bypass this difficulty and leads to a simple and transparent proof.
\end{rem*}

\upskip
\sub{Contractibility of $X_K\lla$ and $X\lla$.}
Now we are ready to fulfill the main goal of this section.

\begin{prop}\label{prop:gRllK}
	The space $X_K\lla$ is contractible for every $\lambda\geq 0$ and every space $X$ from the list \eqref{eq:XgR}.
\end{prop}

\proof
Choose $B\in X_K\lla$.
Then the homotopy $H_t$ defined by \eqref{eq:RKtoB} preserves $X_K\lla$ 
and defines a contraction of $X_K\lla$ to the singleton $\set{B}$. % for every $\lambda\geq 0$.
Indeed,
\[ t\inv u_t Bu_t^* \in X_K(\Ran p_t)\lla \txt{and} (1-t)\inv v_t Av_t^* \in X_K(\Ran q_t)\lla \] 
for every $A,B\in X_K\lla$ and $t\in(0,1)$, 
so $H_t(A)\in X_K\lla$.
\endproof

\begin{prop}\label{prop:gRll}
	The space $X\lla$ is contractible for every $\lambda\geq 0$ and every space $X$ from the list \eqref{eq:XgR}.
\end{prop}

\proof
Choose a positive injective compact operator $k$ of norm $\leq 1$,
and let $H'$ be the corresponding homotopy given by formula \eqref{eq:RtoRK}. 
We will show that $H'_t$ preserves $X\lla$. 
Then the restriction of $H'$ to $X\lla$ provides a deformation of this space to a subspace of $X_K\lla$,
which is contractible by Proposition \ref{prop:gRllK}.

First note that $\norm{c_t \alpha c_t}\leq\norm{\alpha}<\lambda\inv$ 
for every $\alpha$ of norm less than $\lambda\inv$,
and thus $H'_t(A)\in\Rsa\lla$ for every $A\in\Rsa\lla$.
This proves the proposition for $X=\gR\sa$.

An operator $A\in\Rinv\sa$ is essentially positive if any of the following equivalent conditions holds
(as above, we denote $\alpha=A\inv$):
\begin{enumerate}[topsep=-2pt, itemsep=0pt, parsep=3pt, partopsep=0pt]
	\item[(1)] The range $V_-$ of $p=\onem(A)=\onem(\alpha)$ has finite dimension.
	\item[(2)] $\alpha\geq 0$ on some (not necessarily invariant) subspace $V$ of $H$ of finite codimension.
	\item[(3)] There is a finite rank operator $b$ such that $\alpha+b\geq 0$.
\end{enumerate}
Indeed, (1) is in fact the definition of essential positivity, 
for $(1\Rightarrow 3)$ one can take $b = \norm{\alpha}\cdot p$,
for $(3\Rightarrow 2)$ one can take $V = \ker b$,
and $(2\Rightarrow 1)$ follows from $V\cap V_-=0$. %triviality of the intersection $V\cap V_-$,

Suppose that $A\in\Rp\inve$. 
Then $\alpha=A\inv$ satisfies (3), so $h'_t(\alpha) = c_t \alpha c_t$ 
is the difference of a positive operator $c_t(\alpha+b)c_t$ and a finite rank operator $c_t b c_t$
and thus also satisfies (3).
Therefore, $H'_t$ preserves $\Rp\inve$. 
A similar reasoning shows that $H'_t$ preserves $\Rm\inve$.

It remains to show that $H'_t$ preserves $\Rinv\st$.
Let $A\in\Rinv$ and $\alpha=A\inv$. 
Denote by $W_-$ and $W_+$ the ranges of the spectral projections 
$\onem(c_t \alpha c_t)$ and $\one_{(0,\infty)}(c_t \alpha c_t)$. % respectively. 
The sum $c_t(W_-)+c_t(W_+) = c_t(H)$ is dense in $H$.
Suppose that $H'_t(A)\in\Rp\inve$, that is, $W_-$ is finite-dimensional. 
Then $c_t(W_-)$  is also finite-dimensional and hence the closure $V$ of $c_t(W_+)$ has finite codimension.
The operator $\alpha$ is positive on $c_t(W_+)$ and thus also on $V$, 
which implies $A\in\Rp\inve$.
Similarly, $H'_t(A)\in\Rm\inve$ implies $A\in\Rm\inve$.
Therefore, $A\in\Rinv\st = \Rinv\setminus\br{\Rp\inve\cup\Rm\inve}$ implies $H'_t(A)\in\Rinv\st$. 
This completes the proof of the proposition.
\endproof

\section{Proof of Theorem \ref{thm:inv}: odd case}\label{sec:D1}

In this section we prove the second part of Theorem \ref{thm:inv}, 
as well as its analogue for bounded operators which we will need in the proof of Theorem \ref{thm:K0}.

\sub{Passing to odd operators.}
We will use the standard embedding 
\begin{equation}\label{eq:iotaR}
	\iota_{\Reg}\colon\Reg(H)\hookto\Reg\sa(\Hhat), \quad 
	A\mapsto\Ahat = \smatr{0 & A^* \\ A & 0}.
\end{equation}
The map $A\mapsto A^*$ is both Riesz and graph continuous,
since $\f(A^*)=\f(A)^*$ and $\p(A^*)=1-I\cdot\p(A)\cdot I\inv$.
Therefore, \eqref{eq:iotaR} provides homeomorphisms 
\[ \rR(H)\to\rR^1(\Hhat) \txt{and} \gR(H)\to\gR^1(\Hhat), \]
where $\Reg^1(\Hhat)$ denotes the subset of $\Rsa(\Hhat)$ consisting of odd operators
(with respect to the grading symmetry $J = \smatr{1 & 0 \\ 0 & -1}$ of $\Hhat$),
\[ \Reg^1(\Hhat) = \sett{A\in\Rsa(\Hhat)}{JA+AJ=0}. \]
An operator $A$ is invertible (resp. has compact resolvent) if and only if 
$\Ahat$ is invertible (resp. has compact resolvent).
The spaces $X\lla$ and $X_K\lla$ in the statement of Theorem  \ref{thm:inv}(2)
are identified via \eqref{eq:iotaR} with $X^1\lla$ and $X^1_K\lla$ in our current notations.

\begin{prop}\label{prop:gR1}
 The spaces $\gR\lla\cong\gR^1\lla$ and $\gR_K\lla\cong\gR^1_K\lla$ are contractible for every $\lambda\geq 0$.
\end{prop}

\proof
We use homotopies constructed in Section \ref{sec:Dg}, but choose their parameters appropriately.
In the proofs of Propositions \ref{prop:RKinv} and \ref{prop:gRllK}
we replace $u_t$ by $u_t\oplus u_t$, $v_t$ by $v_t\oplus v_t$, and choose $B\in\gR^1_K\lla$.
In the proofs of Propositions \ref{prop:Rinv} and \ref{prop:gRll} we replace $k$ by $k\oplus k$.
Then these homotopies preserve the subspaces of odd operators 
and thus provide contractions of $\gR_K^1\lla$ and $\gR^1\lla$.
\endproof

\sub{Polar decomposition.}
Let $\Reg\gla$ (resp. $\B\gla$) denote the set of regular (resp. bounded) self-adjoint operators in $H$
whose spectrum is contained in the open ray $(\lambda,\infty)$, as in Proposition \ref{prop:conv+}.
The map
\begin{equation}\label{eq:URR+}
	\U(H)\times\Reg^{>0}\to\Rinv, \quad (u,B)\mapsto uB 
\end{equation}
is bijective, with the inverse map given by the polar decomposition of invertible operators:
\[ A\mapsto(\uu_A,|A|), \txt{where} \uu_A\in\U(H) \]
is the unitary part of the polar decomposition.
The restriction of \eqref{eq:URR+} to bounded operators, $\U\times\B^{>0}\to\B\inve$, is a homeomorphism. 
Since $\f(uB)=u\cdot\f(B)$, \eqref{eq:URR+} is a homeomorphism with respect to the Riesz topology on $\Reg$.

By the definition, $\Reg\lla$ consists of regular operators $A$ such that $\Ahat\in\Rsa\lla$.
The last condition is equivalent to $\Ahat^2 = A^*A\oplus AA^*\in\Reg^{>\lambda^2}$,
so $A\in\Reg\lla$ implies $|A|\in\Reg\gla$.
Conversely, $uB\in\Reg\lla$ for every $B\in\Reg\gla$ and $u\in\U$.
Therefore, \eqref{eq:URR+} provides a homeomorphism
\begin{equation}\label{eq:URR>lambda}
	\U(H)\times\rR\gla\to\rR\lla. 
\end{equation}
Restricting to operators with compact resolvents, we obtain a homeomorphism 
\begin{equation*}
	\U(H)\times\rR_K\gla\to\rR_K\lla. 
\end{equation*}
For bounded operators the situation is the same. 
The map $A\mapsto\Ahat$ provides homeomorphisms $\B(H)\to\B^1(\Hhat)$ and $\D(H)\to\D^1(\Hhat)$,
where the superscript $1$ stays for the subspaces of odd self-adjoint operators,
and we denote
\begin{align*}
  \B\lla &= \sett{A\in\B}{\Spec(\Ahat)\cap\lla = \emptyset} \cong\B^1\lla, \\
	\D\lla &= \sett{A\in\D}{\Spec(\Ahat)\cap\f\br{\lla} = \emptyset}\cong\D^1\lla. 
\end{align*}
The restrictions of \eqref{eq:URR+} provide homeomorphisms
\[ \U(H)\times\B\gla\to\B\lla, \quad \U(H)\times\D^{>\mu}\to\D\lla, \txt{and} \U(H)\times\D_K^{>\mu}\to\D_K\lla, \]
where $\mu=\f(\lambda)$.

The unitary group $\U(H)$ is contractible by Kuiper theorem. 
By Proposition \ref{prop:conv+}, $\rR\gla$ and $\rR_K\gla$ are contractible.
A similar reasoning shows that $\B\gla$, $\D^{>\mu}$, and $\D_K^{>\mu}$ are convex and thus contractible.
Hence we obtain the following result.

\begin{prop}\label{prop:RBD1}
 Let $X$ be one of the spaces $\B$, $\rR$, $\rR_K$, $\D$, or $\D_K$.
 Then $X\lla$ is contractible for every $\lambda\geq 0$.
\end{prop}

\upskip
\sub{Unitary operators.}
As we will see in Section \ref{sec:proofA},
for the unitary group the role of the subspace of odd operators plays the space
\begin{equation}\label{eq:U1}
	\U^1 = \sett{u\in\U(\Hhat)}{JuJ=u^*}.
\end{equation}

\begin{prop}\label{prop:U1}
  The subspaces
 $\U^1\lla = \U^1\cap\U\lla$ and $\U^1_K\lla = \U^1\cap\U_K\lla$ of $\U(\Hhat)$
 are both contractible for every $\lambda\geq 0$.
\end{prop}

\proof
The contraction of $\U\lla$ from the proof of Proposition \ref{prop:BDUll}
preserves the subspaces $\U^1\lla$ and $\U^1_K\lla$, so they are also contractible.
\endproof

\section{Proof of Theorems \ref{thm:ess-pos} and \ref{thm:K1}}\label{sec:proofBC}

This section is devoted to the proof of Theorem \ref{thm:ess-pos} and 
the following stronger version of Theorem \ref{thm:K1} (it contains Theorem \ref{thm:K1} as a part).

\begin{thmC'}
	All the maps on Diagram \eqref{diag:K1all} below are homotopy equivalences.
	Consequently, all the spaces on the diagram are classifying spaces for the functor $K^1$.
\begin{equation}\label{diag:K1all}
	\begin{tikzcd}[row sep=small]
		&& \rR\st_K \arrow{rr} \hookr{dd}[near start]{\f} \hookrb{ld} && 
		  \gR\sa_K \hookr{dd}[near start]{\kap} \hookrb{ld}  \\
		\B\st_F \hookr{r} & \rR\st_F \arrow[crossing over]{rr} \hookr{dd}[near start]{\f} &&  \gR\sa_F   \\
		&& \Deu\st \arrow{rr}[near start]{\kat} \hookrb{ld} && 		  \U_K \hookrb{ld}  \\
		& \D\st_F \arrow{rr}{\kat} && 	\U_F \arrow[from=uu, hookrightarrow, crossing over, "{\kap}"{near start}]
	\end{tikzcd}
\end{equation}
\end{thmC'}

\upskip
Here $\kat\colon\D\sa(H)\to\U(H)$ is the map defined by formula \eqref{eq:kat}.
As was explained in Section \ref{sec:maps}, it takes $\Deu\sa$ to $\UK$ and $\D\sa_F$ to $\UF$.

\sub{Plan of the proof.}
The space $\UK(H)$ is well known to be a classifying space for the functor $K^1$.
Therefore, the last statement of Theorem \ref{thm:K1}' follows from the first one.

The proof of Theorem \ref{thm:ess-pos} and the first part of Theorem \ref{thm:K1}' follows the idea described in the Introduction. 
We handle all the maps on Diagrams \eqref{diag:+-*} and \eqref{diag:K1all} at once.

Let $\T$ be the set of all symmetric (with respect to zero) finite non-empty subsets of $\R$. 
The finite union of elements of $\T$ is again an element of $\T$.
For every space $X$ on Diagrams \eqref{diag:+-*} and \eqref{diag:K1all}
we will define an open covering $(X_{\tau})_{\tau\in\T}$ of $X$ indexed by $\T$ 
so that $X_{\tau}\cap X_{\tau'} = X_{\tau\cup\tau'}$
and $\varphi$ takes $X_{\tau}$ to $Y_{\tau}$ for every arrow $\varphi\colon X\to Y$ on the diagrams.
We will show in Proposition \ref{prop:phi-tau-he} 
that the restriction $\varphi_{\tau}\colon X_{\tau}\to Y_{\tau}$ is a homotopy equivalence.
The topology of every such space $X$ is induced by its inclusion to the metric space $\B(H)$,
so $X$ itself is metric and thus every open covering of $X$ is numerable.
Theorems \ref{thm:ess-pos} and \ref{thm:K1}' 
follow immediately from these properties and a theorem of tom Dieck \cite[Theorem 1]{tD}.

\sub{Construction of coverings.}
First notice that every space on Diagrams \eqref{diag:+-*} and \eqref{diag:K1all} 
arises as a subspace of one of the spaces on Diagram \eqref{diag:K1}.
The spaces on Diagram \eqref{diag:K1} are divided naturally into the three groups: 
\begin{enumerate}[topsep=-4pt, itemsep=0pt, parsep=3pt, partopsep=0pt]
	\item[(R)] $\Bsa$, $\rR\sa$ and $\gR\sa$ on the top.
	\item[(D)] The unit ball $\D\sa$ on the bottom left.
	\item[(U)] The unitary group $\U$ on the bottom right.
\end{enumerate}
This determines a partition of the spaces on Diagrams \eqref{diag:+-*} and \eqref{diag:K1all} into the three groups:
\begin{enumerate}[topsep=-4pt, itemsep=0pt, parsep=3pt, partopsep=0pt]
	\item[(R)] Subspaces of $\Bsa$, $\rR\sa$ and $\gR\sa$; 
these are the spaces on the top part of \eqref{diag:K1all} and all the spaces of \eqref{diag:+-*}.
	\item[(D)] Subspaces of $\D\sa$; these are the spaces on the left bottom part of \eqref{diag:K1all}.
	\item[(U)] Subspaces of $\U$; these are the spaces on the right bottom part of \eqref{diag:K1all}.
\end{enumerate}
The arrows inside each group are natural inclusions.
The arrow from a space of type (R) to a space of type (D) is given by the bounded transfrom $\f$,
from (R) to (U) by the Cayley transform $\kap$, 
and from (D) to (U) by $\kat$.

Let $\taub$ denote the convex hull of $\tau\in\T$; it is a closed symmetric interval in $\R$.
In order to simplify notations, we will use the symbols $\sigma$ and $\sigb$ as replacements for:
\begin{itemize}[topsep=-4pt, itemsep=0pt, parsep=3pt, partopsep=0pt]
	\item $\tau$ and $\taub$ if we deal with a space of type (R),
	\item $\f(\tau)$ and $\f(\taub)$ if we deal with a space of type (D),
	\item $\kap(\tau)$ and $\kap(\taub)$ if we deal with a space of type (U).
\end{itemize}
We will consider the spaces on the diagrams as depending on a Hilbert space $H$
and will write $X(H)$ when the argument $H$ is not fixed. 

If $X$ is one of the spaces on Diagram \eqref{diag:+-*} or \eqref{diag:K1all},
then we define $X_{\tau}$ by the formula
\begin{equation}\label{eq:Xtau}
	X_{\tau} = \sett{A\in X}{\Spec(A)\cap\sigma = \emptyset = \Spec\ess(A)\cap\sigb},
\end{equation}
where $\Spec\ess(A)$ is the essential spectrum of $A$.

\begin{rem*}
Of course, if $X$ is a subspace of $\Rsa_K$, $\Deu\sa$, or $\UK$, 
then the condition for the essential spectrum in \eqref{eq:Xtau} is void. 
We prefer to write the definition of $X_{\tau}$ in the same form for all possible spaces $X$, though, 
because we will work with all these spaces simultaneously.
\end{rem*}

\upskip
\sub{Properties of coverings.}
Clearly, $\cup X_{\tau} = X$, $X_{\tau}\cap X_{\tau'} = X_{\tau\cup\tau'}$,
and $\varphi(X_{\tau})\subset Y_{\tau}$ for every arrow $\varphi\colon X\to Y$ 
on the diagrams.

\begin{prop}\label{prop:open}
	Each set $X_{\tau}$ is open in $X$ and the map $\one_{\sigb}\colon X_{\tau}\to\Proj(H)$ is continuous.
\end{prop}

\proof
The proof is a standard application of perturbation results from \cite{Kato}.
We write it for the case $X\subset\gR\sa$;
for bounded operators the proof is completely similar (with $\tau$ replaced by $\sigma$ 
and, in the unitary case, $\R$ replaced by the unit circle in $\CC$),
and for $X\subset\rR\sa$ one needs to apply the bounded transform first.

If $\tau$ is a singleton, $\tau=\set{0}$, then $X_{\tau}$ consists of invertible operators 
and is open in $X$ by \cite[Theorem IV.3.1]{Kato}; in this case $\one_{\sigb}\equiv 0$ on $X_{\tau}$.

Let $\tau = \set{\tau_1<\ldots<\tau_n}$.
Choose a compact region $G$ in $\CC$ bounded by a simple closed smooth curve $\Gamma$ such that 
$G\cap\R=\brr{\tau_1,\tau_n}$ and $\Gamma\cap\R=\set{\tau_1,\tau_n}$.
Let $A_0$ be an arbitrary element of $X_{\tau}$. 
By \cite[Theorem IV.3.1]{Kato}, there is an open neighborhood $U$ of $A_0$ in $\gR\sa$ such that 
\[ \Spec(A)\cap\tau = \emptyset \txt{for every} A\in U. \]
In particular, $\Spec(A)\cap\Gamma = \emptyset$.
The curve $\Gamma$ separates the parts of the spectrum of $A$ lying inside and outside of $\Gamma$ 
and determines the corresponding %(orthogonal) 
decomposition $H=H'_A\oplus H''_A$ into invariant subspaces of $A$ \cite[Theorem III.6.17]{Kato}.
Since $A$ is self-adjoint, $\Spec(A)\cap G = \Spec(A)\cap\taub$, 
so the invariant subspace $H'_A$ is the range of the spectral projection 
$\one_{\taub}(A) = \one_{\sigb}(A)$.
By \cite[Theorem IV.3.16]{Kato}, this projection depends continuously on $A\in U$.
Decreasing $U$ if necessary, we get $\dim H'_A=\dim H'_{A_0}<\infty$ for $A\in U$,
so the second condition in \eqref{eq:Xtau} is also satisfied for such $A$. 
Therefore, $U\subset X_\tau$ and thus $X_\tau$ is open.
\endproof

To fulfill the program described in the beginning of the section, it remains to prove the following statement:

\begin{prop}\label{prop:phi-tau-he}
For every $\tau\in\T$ and every arrow $\varphi\colon X\to Y$ on Diagrams \eqref{diag:+-*} and \eqref{diag:K1all},
the restriction $\varphi_{\tau}\colon X_{\tau}\to Y_{\tau}$ is a homotopy equivalence. 
\end{prop}

\upskip
The rest of the section is devoted to the proof of this proposition.

\sub{Fiber bundle structure.}
Let $X$ be one of the spaces on the diagrams and $\tau\in\T$.
By Proposition \ref{prop:open}, the spectral projection 
$\one_{\sigb}\colon X_{\tau}\to\Proj(H)$ is continuous.
The image of this map is the subspace $\Pfin$ of $\Proj(H)$ consisting of projections of finite rank, so we obtain a continuous map
\begin{equation}\label{eq:Xtau-P}
 \one_{\sigb}\colon X_{\tau}\to\Pfin. 
\end{equation}

\begin{lem}\label{lem:Xtau-trivial}
	The map \eqref{eq:Xtau-P} is a locally trivial fiber bundle.
\end{lem}

\proof
For an arbitrary $q\in\Pfin$, choose a continuous map $g$ from a neighborhood $W\subset\Pfin$ of $q$ to $\U(H)$
such that $g(p)\cdot p\cdot g(p)^* \equiv q$ for every $p\in W$, as in \eqref{eq:gWU}.
Let $\restr{X_{\tau}}{W}$ and $\restr{X_{\tau}}{q}$ denote the inverse images under \eqref{eq:Xtau-P} 
of $W$ and $q$ respectively.
The unitary group $\U(H)$ acts continuously on both $X_{\tau}$ and $\Pfin$ by conjugations
and \eqref{eq:Xtau-P} is equivariant with respect to this action.
Hence the map
\begin{equation*}
	\restr{X_{\tau}}{W} \to \restr{X_{\tau}}{q}\times W, \quad 
		A\mapsto \br{ g(p) A g(p)^*,\, p }, \txt{where} p=\one_{\sigb}(A),
\end{equation*}
is a bundle isomorphism providing a trivialization of \eqref{eq:Xtau-P} over $W$.
\endproof

Let $H\times\Pfin = \H'\oplus\H''$ be the restriction of the canonical decomposition \eqref{eq:HHH} to $\Pfin$.
We define $X'_{\tau}$ and $X''_{\tau}$ to be the (locally trivial) fiber bundles over $\Pfin$ 
associated with $\H'$ and $\H''$, whose fibers over $p\in\Pfin$ are given by the formulas
\begin{equation*}\label{eq:X'X''}
		X'_{\tau,\,p} = \sett{A\in X(\H'_p)}{\Spec(A)\subset\sigb\setminus\sigma} \;\; \text{and} \;\;
		X''_{\tau,\,p} = \sett{A\in X(\H''_p)}{\Spec(A)\cap\sigb = \emptyset}.
\end{equation*}

\begin{lem}\label{lem:Xtau}
The bundle $X_{\tau}\to\Pfin$ is naturally decomposed as the fiber product
\begin{equation}\label{eq:XXXP}
  X_{\tau} = X'_{\tau}\times_{\Pfin}X''_{\tau},
\end{equation}
where the bundle maps $X_{\tau}\to X'_{\tau}$ and $X_{\tau}\to X''_{\tau}$ over $\Pfin$
are given by the restriction of $A\in X_{\tau}$ to the range and the kernel of $\one_{\sigb}(A)$ respectively.
The product \eqref{eq:XXXP} is functorial by $X$ in the following sense: 
every arrow $\varphi\colon X\to Y$ on Diagrams \eqref{diag:+-*} and \eqref{diag:K1all}
induces a commutative diagram of fiber bundles over $\Pfin$:
\begin{equation}\label{eq:XXXP2}
	\begin{tikzcd}
		X'_{\tau} \arrow{d}{\varphi'_{\tau}} & 
		  X_{\tau} \arrow{l} \arrow{r} \arrow{d}{\varphi_{\tau}} & 
			  X''_{\tau} \arrow{d}{\varphi''_{\tau}} \\
		Y'_{\taub} & Y_{\tau} \arrow{l} \arrow{r} & Y''_{\tau} 
 \end{tikzcd}
\end{equation}
The induced map $\varphi'_{\tau}\colon X'_{\tau}\to Y'_{\tau}$ is a bundle isomorphism over $\Pfin$.
\end{lem}

\proof
By the definition of $\sigb$, the spectral projections $\one_{\sigb}(\varphi(A))$ and $\one_{\sigb}(A)$ coincide,
so \eqref{eq:XXXP2} is a commutative diagram of fiber bundles.
Since all the bundles in \eqref{eq:XXXP} and \eqref{eq:XXXP2} are locally trivial, 
we only need to check these statements fiberwise, over each $p\in\Pfin$. 
But fiberwise the first two statements are obvious.
For the last statement, note that 
\begin{equation*}
  X'_{\tau,\,p} = \case{
	     \sett{A\in \Bsa(\H'_p)}{\Spec(A)\subset\taub\setminus\tau} &\text{ for } X \text{ of type (R)}, \\
			 \sett{A\in \D\sa(\H'_p)}{\Spec(A)\subset\f(\taub\setminus\tau)} &\text{ for } X \text{ of type (D)}, \\
			\sett{A\in \U(\H'_p)}{\Spec(A)\subset\kap(\taub\setminus\tau)} &\text{ for } X \text{ of type (U)}
			}
\end{equation*}
If $X$ and $Y$ are of the same type, then $X'_{\tau}$ and $Y'_{\tau}$ coincide and $\varphi'_{\tau}$ is the identity.
If $X$ is of type (R) and $Y$ is of type (D), 
then $\varphi = \f$ and the map $\f\colon X'_{\tau,\,p}\to Y'_{\tau,\,p}$ is obviously a homeomorphism.
The other two cases, $\varphi = \kap$ and $\varphi = \kat$, are similar.
\endproof

\begin{lem}\label{lem:X''tr}
  The fiber bundle $X''_{\tau}\to\Pfin$ is trivial 
	(that is, isomorphic to the trivial bundle $X''_{\tau,\,0}\times \Pfin\to\Pfin$)
	for every space $X$ on Diagrams \eqref{diag:+-*} and \eqref{diag:K1all} and every $\tau\in\T$.
\end{lem}

\proof
The locally trivial Hilbert bundle $\H''$ over the paracompact space $\Pfin$ 
has infinite-dimensional separable fibers and thus is trivial by Kuiper's theorem.
Therefore, the fiber bundle $X''_{\tau}\to\Pfin$ associated with $\H''$ is also trivial.
\endproof

\sub{Proof of Proposition \ref{prop:phi-tau-he}.}
Lemmas \ref{lem:Xtau} and \ref{lem:X''tr} imply that the map $\varphi_{\tau}\colon X_{\tau}\to Y_{\tau}$
is the product of $\varphi'_{\tau}\colon X'_{\tau}\to Y'_{\tau}$ and $\varphi''_{\tau}\colon X''_{\tau,\,0}\to Y''_{\tau,\,0}$.
The first factor $\varphi'_{\tau}$ is a homeomorphism. 
The spaces $X''_{\tau,\,0} = X\lla$ and $Y''_{\tau,\,0} = Y\lla$, where $\lla=\taub$, 
are contractible by Propositions \ref{prop:rRll}, \ref{prop:BDUll}, \ref{prop:gRllK}, and \ref{prop:gRll},
and thus the second factor $X''_{\tau,\,0}\to Y''_{\tau,\,0}$ is a homotopy equivalence.
Therefore, the product $\varphi_{\tau}\colon X_{\tau}\to Y_{\tau}$ is also a homotopy equivalence.
This completes the proof of the proposition and Theorems \ref{thm:ess-pos} and \ref{thm:K1}'.
\endproof

\section{Proof of Theorem \ref{thm:K0}}\label{sec:proofA}

This section is devoted to the proof of the following theorem, which contains Theorem \ref{thm:K0} as a part. 

\begin{thmA'}
	All the maps on the Diagram \eqref{diag:K0all} below are homotopy equivalences.
	Consequently, all the spaces on the diagram are classifying spaces for the functor $K^0$.
\begin{equation}\label{diag:K0all}
	\begin{tikzcd}[row sep=small]
		&& \rR_K \arrow{rr} \hookr{dd}[near start]{\f} \hookrb{ld} && 
		  \gR_K \hookr{dd}[near start]{\p} \hookrb{ld}  \\
		\B_F \hookr{r} & \rR_F \arrow[crossing over]{rr} \hookr{dd}[near start]{\f} &&  \gR_F   \\
		&& \Deu \arrow{rr}[near start]{\pt} \hookrb{ld} && 		  \PK \hookrb{ld}  \\
		& \D_F \arrow{rr}{\pt} && 	\PF \arrow[from=uu, hookrightarrow, crossing over, "{\p}"{near start}]
	\end{tikzcd}
\end{equation}
\end{thmA'}

\upskip
Here $\pt\colon\D(H)\to\Proj(\Hhat)$ is the map defined by formula \eqref{eq:pt}.
As was explained in Section \ref{sec:maps}, it takes $\Deu$ to $\PK$ and $\D_F$ to $\PF$.

\medskip
\proof
Let us describe briefly the plan of the proof.
Let $D'$ denote Diagram \eqref{diag:K0all} and $D$ denote Diagram \eqref{diag:K1all}.
Consider the diagram $D(\Hhat)$ corresponding to the $\Z/2$-graded Hilbert space $\Hhat = H\oplus H$.
There is a natural morphism of the diagrams $\iota\colon D'(H)\to D(\Hhat)$, 
which induces an embedding $\iota\colon X'\hookto X$ of each space $X'$ of $D'(H)$ to the corresponding space $X$ of $D(\Hhat)$.
This embedding takes $X'$ to the subspace $X^1 = \iota(X')$ of $X$.
Theorem \ref{thm:K0}' is equivalent to the following statement: 
for every arrow $\varphi\colon X\to Y$ of $D(\Hhat)$, its restriction $\varphi^1\colon X^1\to Y^1$ is a homotopy equivalence.
By \cite[Theorem 1]{tD}, it is sufficient to show that 
$X^1_{\tau} = X^1\cap X_{\tau} \to Y^1_{\tau} = Y^1\cap Y_{\tau}$ 
is a homotopy equivalence for every $\tau\in\T$.
The proof of this mostly follows the proof of Theorem \ref{thm:K1}', 
and we keep designations of Sections \ref{sec:D1} and \ref{sec:proofBC}.
Once we show that all the maps are homotopy equivalences, the last statement of the theorem 
follows from a result of Atiyah and J\"anich about the homotopy type of $\BF(H)$ \cite[Theorem A1]{Atiyah}.

Let us describe this construction in more detail.
We will use the standard embedding \eqref{eq:iotaR}, 
which provides a homeomorphism (both in the Riesz and graph topology on these two spaces)
of $\Reg(H)$ onto the subspace $\Reg^1(\Hhat)$ of $\Reg\sa(\Hhat)$ consisting of odd operators.
If $A\in\RF$, then $\Ahat\in\RF\st$.
The bounded transform commutes with $\iota_{\Reg}$.
The same map $A\mapsto\Ahat$ takes $\B(H)$ onto $\B^1(\Hhat)$ and $\D(H)$ onto $\D^1(\Hhat)$.
We denote the corresponding embedding $\D(H)\hookto\D\sa(\Hhat)$ by $\iota_{\D}$.

The restriction of $\kat\colon\D\sa(\Hhat)\to\U(\Hhat)$ to $\D^1(\Hhat)$ is given by the formula
\begin{equation*}
	\kat(\hat{a}) = \matr{2a^*a-1 & -2ia^*\sqrt{1-aa^*} \\ -2ia\sqrt{1-a^*a} & 2aa^*-1}
		= v\cdot\br{2\pt(a)-1}\cdot v \txt{for} a\in\D(H), 
\end{equation*}
where $v = \smatr{i & 0 \\ 0 & -1}\in\U(\Hhat)$, $v^2=-J$.
For an arbitrary $r\in\B(\Hhat)$ and $u = vrv$ we have % = J(v\inv rv)
\[ Ju = Jvrv = -v\inv rv \txt{and} (Ju)^2 = v\inv r^2\,v. \]
Hence $r^2=1$ if and only if $(Ju)^2 = 1$. 
Clearly, $r$ is unitary if and only if $u$ is unitary.
Therefore, we obtain a homeomorphism
\begin{equation}
	\iota_{\Proj}\colon \Proj(\Hhat)\to\U^1(\Hhat) = \sett{u\in\U(\Hhat)}{(Ju)^2=1}, 
	\quad \iota_{\Proj}(p) = v(2p-1)v,
\end{equation}
which makes the following diagram commutative:
\begin{equation*}
	\begin{tikzcd}
		\D(H) \arrow{r}{\pt} \arrow{d}{\iota_{\D}} &
		  \Proj(\Hhat) \arrow{d}{\iota_{\Proj}} \\
		\D^1(\Hhat) \arrow{r}{\kat} &
		  \U^1(\Hhat)
	\end{tikzcd}
\end{equation*}
Since $\iota_{\Proj}(\pinf)=1$ and $\iota_{\Proj}(p_0)=-1$, 
the homeomorphism $\iota_{\Proj}$ takes $\PK$ to $\UK^1=\UK\cap\U^1$ and $\PF$ to $\UF^1=\UF\cap\U^1$.

Let $X$ be a space from the diagram $D(\Hhat)$ and $\tau\in\T$.
Since the interval $\taub$ is symmetric with respect to zero, 
the restriction of \eqref{eq:Xtau-P} to odd operators provides the map
\begin{equation}\label{eq:X1tau}
 \one_{\sigb}\colon X^1_{\tau}\to\Pfin^0 = \sett{p\in\Pfin(\Hhat)}{pJ=Jp}. 
\end{equation}

\begin{lem}\label{lem:X1tau}
	The map \eqref{eq:X1tau} is a locally trivial fiber bundle.
\end{lem}

\proof
The subgroup $\U^0(\Hhat)=\sett{u\in\U(\Hhat)}{uJ=Ju}$ of $\U(\Hhat)$ 
acts continuously on both $X^1_{\tau}$ and $\Pfin^0$ by conjugation,
and \eqref{eq:X1tau} is equivariant with respect to this action.
Every element $q$ of $\Pfin^0$ has the form $q=q'\oplus q''$, where $q',q''\in\Pfin(H)$. 
Let $g'\colon W'\to\U(H)$ and $g''\colon W''\to\U(H)$ be maps
from neighborhoods $W'\ni q'$ and $W''\ni q''$, as in \eqref{eq:gWU}.
Then $W=W'\times W''$ is a neighborhood of $q$ in $\Pfin^0$ and the map
\[ g\colon W\to\U^0(\Hhat), \quad g(p'\oplus p'')=g'(p')\oplus g''(p'') \]
satisfies the condition $g(p)pg(p)^*\equiv q$ for every $p\in W$.
Therefore, the map 
\begin{equation*}
	\restr{X^1_{\tau}}{W} \to \restr{X^1_{\tau}}{q}\times W, \quad 
		A\mapsto \br{ g(p) A g(p)^*,\, p }, \txt{where} p=\one_{\sigb}(A).
\end{equation*}
provides a trivialization of \eqref{eq:X1tau} over $W$, 
similarly to the proof of Lemma \ref{lem:Xtau-trivial}.
\endproof

The rest of the proof of Theorem \ref{thm:K0}'
is completely similar to the proof of Theorem \ref{thm:K1}',
but the grading of $\Hhat$ should be taken into account.
The fiber bundle \eqref{eq:X1tau} is decomposed into the fiber product, as in \eqref{eq:XXXP}.
The map $\varphi^1\colon X^1_{\tau}\to Y^1_{\tau}$ respects these bundle and fiber product structures 
and induces a bundle isomorphism on the first factor of the fiber product.
The second factor is a trivial fiber bundle over $\Pfin^0$ with the fiber $X^1\lla = X^1\cap X\lla$.
It only remains to apply Propositions \ref{prop:gR1}, \ref{prop:RBD1}, and \ref{prop:U1},
where contractibility of $X^1\lla$ is proven for every $\lambda\geq 0$ 
and every space $X$ of the diagram $D(\Hhat)$.
\endproof

\section{Essentially positive / negative operators}\label{sec:ess-pos}

\upskip
\begin{thm}\label{thm:dens}
Each of the subspaces $\RKm(H)$, $\RKp(H)$, and $\RK\st(H)$ is dense in $\RK\sa(H)$ in the graph topology.
\end{thm}

\proof
Let $A\in\RK\sa(H)$ and $\eps>0$.
Choose $c>0$ large enough, so that $|\kap(c)-1|<\eps/2$.
The projection $p=\one_{[-c,\,c]}(A)$ has finite rank, 
so the kernel $H'$ of $p$ is infinite-dimensional.
Suppose that the resolvent set of an operator $B\in\RK\sa(H')$ contains $[-c,c]$.
Then $A' = Ap+B(1-p)$ is an element of $\RK\sa(H)$ and $\norm{\kap(A')-\kap(A)}<\eps$.
Moreover, $A'$ coincides with $B$ on the invariant subspace $H'$ of finite codimension.
Choosing such a $B$ from $\RKm(H')$, $\RKp(H')$, or $\RK\st(H')$,
we get $A'\in\RKm(H)$, $\RKp(H)$, or $\RK\st(H)$ respectively.
Since $A$ and $\eps$ were chosen arbitrarily, this completes the proof of the theorem.
\endproof

\begin{rem*}
The proof of Theorem \ref{thm:dens} uses essentially the fact that $A$ has compact resolvent.
This theorem has no analogue for the space $\RF\sa(H)$ of \emph{Fredholm} operators.
This can be explained loosely as follows.
Roughly speaking, a self-adjoint operator with compact resolvent has all its essential spectrum at the infinity
(that is, its Cayley transform has all its essential spectrum at $1$).
This part of the spectrum can be smuggled through the infinity 
from a positive half-neighborhood of the infinity to a negative one and vice versa.
In contrast with that, a Fredholm operator may have a non-empty essential spectrum in a finite part of $\R$.
When such an operator is slightly deformed, this part of the essential spectrum remains near its original position 
and thus cannot be moved through the infinity.
\end{rem*}

\begin{exm}[\textbf{Non-zero spectral flow for essentially positive operators}]\label{exm:sf+}
\label{ex:delta}
Let $A = -d^2/dt^2$ be the Laplace operator acting on complex-valued functions 
$\psi\colon I=[0,1]\to\CC$.
We consider a family of boundary value problems for $A$ 
parametrized by points of the real projective line $\RP^1\cong S^1$.
For $\x=[x_0:x_1]\in\RP^1$, let $\A_{\x}$ be the operator $A$ with the domain 
\begin{equation}\label{eq:delta}
	\dom(\A_{\x}) = \sett{\psi\in H^2(I)}{\psi(0) = 0, \; x_0\psi(1) - x_1\psi'(1) = 0}
\end{equation}
given by local boundary conditions.
Here $\psi' = d\psi/dt$, the type of the boundary condition at $t=1$
is Dirichlet for $\x=[1:0]$, Neumann for $\x=[0:1]$, 
and mixed for all the other values of parameter\footnote{\hspace{1pt}
The family $\A$ restricted to the affine line $\R = \set{[1:x]}\subset\RP^1$ 
coincides with Rellich's example \cite[Example V.4.14]{Kato}. 
It is instructive to look at \cite[Fig. V.4.1]{Kato}.}. 

Every $\A_{\x}$ considered as an unbounded operator on $H = L^2(I)$ 
is a self-adjoint operator with compact resolvent.
As we will see below, all $\A_{\x}$ are essentially positive.

The family $\A=(\A_{\x})$ of regular operators on $H$ is graph continuous.
Indeed, the map $H^2(I)\to\CC^3$, $\psi\mapsto(\psi(0),\psi(1),\psi'(1))$, is continuous,
so $\dom(\A_{\x})$ is a closed subspace of $H^2(I)$
depending continuously on $\x$ in the gap topology on $\Gr\br{H^2(I)}$.
By \cite[Theorem 4.1]{Alonso-Simon}, the map $\A\colon\RP^1\to\RK\sa(H)$ is graph continuous.

Let us look at the spectral graph 
\[ \Gamma=\sett{(\x,\lambda)}{\lambda\in\Spec(\A_{\x})}\in\RP^1\times\R \]
of this operator family.
The space of solutions of the equation $A\psi = \lambda\psi$ 
satisfying the first boundary condition $\psi(0) = 0$ is one-dimensional for every $\lambda\in\R$.
$\Gamma$ intersects the level $\lambda=0$ only at $\x=[1:1]$; 
the corresponding eigenfunction is $\psi(t)=t$.
Similarly, a positive $\lambda=\mu^2>0$ is an eigenvalue of $\A_{\x}$ for exactly one value of $\x$, 
namely $\x = [(e^{2i\mu}+1) : (e^{2i\mu}-1)]$;
the corresponding eigenfunction is $\psi(t) = e^{\mu it} - e^{-\mu it}$. % and has multiplicity 1.
For $\lambda=-\mu^2<0$ the eigenfunction %of $A$ 
is $\psi(t) = e^{\mu t} - e^{-\mu t}$; %, $c=\const$.
it satisfies the second boundary condition if and only if 
\begin{equation}\label{eq:xmu}
	x_0(e^{2\mu}-1) = x_1\mu(e^{2\mu}+1).
\end{equation}
For $\lambda$ running the ray $(-\infty,0)$, $\mu$ runs the positive half of the real axis from $+\infty$ to $0$,
and \eqref{eq:xmu} has exactly one solution $\x(\mu)\in\RP^1$ for every such value $\mu\neq 0$ 
(we can consider only half of the real axes, 
since the values $\mu$ and $-\mu$ give the same eigenspace $\bra{\psi}$ of $\A_{\x}$.)

Combining all this together, we see that the part of $\Gamma$ corresponding to positive values of $\lambda$ 
is an infinite spiral line starting at the point $([1:1], 0)$
and making an infinite number of rotations in the upward direction over the circle $\RP^1$.
The fiber of $\Gamma$ over $\x=[0:1]$ is $\sett{\pi^2 (n+1/2)^2}{n\in\Z, n\geq 0}$,
while the fiber over $\x=[1:0]$ is $\sett{\pi^2 n^2}{n\in\Z, n>0}$.
In contrast with that, the projection of the negative part of $\Gamma$ to $\RP^1$
is the interval $\sett{[1:x]}{x\in(0,1)}$;
the negative eigenvalue $\lambda$ goes to $-\infty$ as $x$ goes to $+0$.

In other words, the graph $\Gamma$ is a single line, 
which has a vertical asymptotics $\lambda\to-\infty$ for $\x$ approaching $[1:+0]$
and the spiral behavior for positive $\lambda$.
The projection to the $\lambda$-axis monotonically increases along the whole line $\Gamma$.
Every operator $\A_{\x}$ is bounded from below and thus is essentially positive;
it has a negative eigenvalue if and only if $\x=[1:x]$ with $x\in(0,1)$.

This description implies that the spectral flow of the graph continuous loop %$\A$ is equal to 1.
$\A\colon\RP^1\to\RK^+(H)$ is equal to 1.
In particular, this loop is not contractible in the graph topology, 
not only in $\gR_K^+(H)$ but also in $\gR_F\sa(H)$.
In contrast with this, every \emph{Riesz} continuous loop of essentially positive  operators is contractible.
It follows that our family $\A$ is not Riesz continuous.

The Riesz discontinuity of $\A$ can be also observed in a more direct way.
The negative eigenvalue of $a_{\x}=\f(\A_{\x})$ goes to $-1$ as $x\to [1:+0]$.
However, $A_{[1:0]}$ and thus $a_{[1:0]}$ is positive.
Therefore, $a=\f\circ\A$ is discontinuous at $\x=[1:0]$ and thus $\A$ is Riesz discontinuous at the same point.
It can be easily seen, for example using \cite[Proposition 11.3]{Pr20}, 
that $\A$ is Riesz continuous at the rest of $\RP^1$.

Let us look more closely at what happens with the eigenfunction when $\x\to[1:+0]$ and $\lambda=-\mu^2\to -\infty$.
The normalized eigenfunction is $\psi(t) = c(e^{\mu t} - e^{-\mu t})$, 
where the constant $c$ is defined by the condition $\norm{\psi}_{L^2}=1$.
When $\mu\to +\infty$, the asymptotics are $c\approx\sqrt{2\mu}\,e^{-\mu}$
and $\psi(t)\approx \sqrt{2\mu}\,e^{\mu(t-1)}$.
Therefore, $\psi$ is more and more concentrating near the right end $t=1$ of the interval $[0,1]$ 
with the increasing of $\mu$.
In the limit, this eigenfunction collapses to the delta-function $\delta(1)$,
which is not an element of $L^2[0,1]$.
\end{exm}

\section*{Conventions and notations}
%\addcontsec{\hspace{18pt}Conventions and notations}

{\small
Throughout the paper, a \q{Hilbert space} always means a separable complex Hilbert space of infinite dimension.

$\B(H)$ is the space of bounded linear operators on a Hilbert space $H$ with the norm topology. 
\\ $\K(H)$ is the subspace of $\B(H)$ consisting of compact operators.
\\ $\Proj(H)$ is the subspace of $\B(H)$ consisting of projections (that is, self-adjoint idempotents). 
\\ $\U(H)$ is the subspace of $\B(H)$ consisting of unitary operators.
\\ $\D(H)$ is the closed unit ball in $\B(H)$.
\\ $\Reg(H)$ is the set of regular (that is, closed and densely defined) operators on $H$.

All the subspaces of $\B(H)$ are supposed to be equipped with the norm topology, if the converse is not stated explicitly.
We use the left superscript for pointing to topology on $\Reg(H)$:
\\ r -- Riesz topology
\\ g -- graph topology 

The right sub- and superscript point to the type of a subspace.
In particular,
\\ $\RF(H)$ is the subspace of $\Reg(H)$ consisting of Fredholm operators.
\\ $\RK(H)$ is the subspace of $\Reg(H)$ consisting of operators with compact resolvents.
\\ $\D_K(H)$ is the subspace of $\D(H)$ consisting of essentially unitary operators.
\\ $\Bsa$, $\Rsa$, $\RK\sa$, ... are the subspaces of self-adjoint operators.
\\ $\Bp$, $\Rp$ are the subspaces of essentially positive self-adjoint operators.
\\ $\Bm$, $\Rm$ are the subspaces of essentially negative self-adjoint operators.
\\ $\B\st$, $\Reg\st$ are the subspaces of self-adjoint operators 
that are neither essentially positive nor essentially negative.
\\ $\B^1(\Hhat)$, $\D^1(\Hhat)$, $\Reg^1(\Hhat)$ are the subspaces of odd self-adjoint operators 
in $\Hhat= H\oplus H$.

}

\end{document}